\documentclass[11pt]{article}   

\usepackage{amsmath, amssymb}    
\usepackage{amsthm}	
\usepackage{amscd} 
\usepackage{graphicx}

\usepackage[utf8]{inputenc}		
\usepackage[T1]{fontenc}		 

\DeclareMathAlphabet{\mathpzc}{OT1}{pzc}{m}{it}

\usepackage{enumitem}

\usepackage{xcolor} 

\usepackage[colorlinks=true,linkcolor=blue,citecolor=magenta]{hyperref}       

 \title{\textbf{The Intermediate Disorder Regime for Brownian Directed Polymers in Poisson Environment}}

\newcommand{\IP}{{\mathbb P}}
\newcommand{\DP}{{\mathrm P}}
\newcommand*\Z{Z_{n}^\omega(\beta)}
\renewcommand{\b}{\beta}
\newcommand*\Hn{H_n^\omega}
\newcommand*\Zc{\mathcal{Z}}
\newcommand*\R{\mathbb{R}}
\newcommand*\Skn{\mathcal{S}_k^n}
\newcommand{\overbar}[1]{\mkern 1.5mu\overline{\mkern-1.5mu#1\mkern-1.5mu}\mkern 1.5mu}
\newcommand{\cvlaw}{\stackrel{{ (d)}}{\longrightarrow}}
\newcommand{\eqlaw}{\stackrel{\rm{ law}}{=}}
\newcommand*\cvLdeux{\overset{L^2}{\longrightarrow}}
\newcommand{\cvto}[2]{\underset{#1 \to #2}{\longrightarrow}}
\newcommand{\om}{\omega}
\newcommand{\rmd}{\mathrm{d}}
\def\P{{\mathbb P}}

\newtheorem{proposition}{Proposition}[section]
\newtheorem{theorem}{Theorem}[section]
\newtheorem{lemma}{Lemma}[section]
\newtheorem{corollary}[theorem]{Corollary}

\newtheorem*{theorem*}{Theorem}
\newtheorem*{lemme*}{Lemme}

\newtheorem*{proposition*}{Proposition}
\newtheorem{conjecture}[theorem]{Conjecture}
\newtheorem{property}[theorem]{Property}

\theoremstyle{definition}
\newtheorem{definition}{Definition}[section]

\theoremstyle{remark}
\newtheorem{remark}{Remark}[section]

\AtBeginDocument{
   \def\MR#1{}  }
   
\makeatletter
\def\blfootnote{\xdef\@thefnmark{}\@footnotetext}
\makeatother

\author{Cl\'ement Cosco\\
   Laboratoire de Probabilit\'es, Statistiques et Mod\'elisation,\\
   Universit\'e Paris Diderot,\\
   \texttt{clement.cosco@lpsm.paris}
}
\date{}

\usepackage[numbers]{natbib}
\begin{document}
\maketitle
\begin{abstract}
We consider the Brownian directed polymer in Poissonian environment in dimension 1+1, under the so-called \emph{intermediate disorder regime} \cite{AlbertsKhaninQuastelAP}, which is a crossover regime between the strong and weak disorder regions. We show that, under a diffusive scaling involving different parameters of the system, the renormalized point-to-point partition function of the polymer converges in law to the solution of the stochastic heat equation with Gaussian multiplicative noise. The Poissonian environment provides a natural setting and strong tools, such as the Wiener-It\^o chaos expansion \cite{LastPenrose}, which, applied to the partition function, is the basic ingredient of the proof.
\blfootnote{\textit{AMS 2010 subject classifications}. 60F05, 82D60.} \blfootnote{\textbf{Key words and phrases}. Directed polymers, Near-critical scaling limits, Chaos expansions, KPZ equation.}
\end{abstract}

\newpage
\tableofcontents

\section{Introduction}

The model of directed polymers in random environment is a simple description of stretched chains which are sensitive to external random impurities. It was first introduced in the physics literature in \cite{HuHe85} and several variations of the model have been studied ever since \cite{ImbrieSpencer88,CY05,
AlbertsKhaninQuastelCRP,OConnellYor}. Although the model has attracted significant attention, many expected results are still far from being proved with mathematical rigor. One can refer to \cite{CStFlour} for a recent review in the discrete setting. In the present paper, we consider a Brownian path directed polymer, where the external impurities are represented by a Poisson point process.

\subsection{The model and its context}

Let $((B_t)_{t\geq 0}, {\DP}_x)$ denote the Brownian motion starting from $x\in\mathbb{R}^d$ and set $\DP=\DP_0$. The \emph{environment} is a Poisson point process $\omega = \sum_i \delta_{s_i,x_i}$ on $[0,\infty)\times\mathbb{R}^d$ of intensity measure $\nu \rmd s \rmd x$. We assume that $\om$ is defined on some probability space $(\Omega,\mathcal{G},\IP)$ and we will denote by $\om_t=\om_{|[0,t]\times\mathbb{R}^d}$ the restriction of the environment to times before $t\geq 0$. Fix $r>0$ and let $U(x)$ be the ball of volume $r^d$, centered at $x\in\mathbb{R}^d$. Define $V_t(B)$ as the tube around path $B$:
\begin{equation}
V_t(B) = \{(s,x) : s\in(0,t],x\in U(B_s)\}.
\end{equation}
For any fixed path $B$, the variable $\om(V_t)$ stands for the energy of the path until time $t$ and corresponds to the number of points that the path encounters; it has the law of a Poisson variable of mean $\nu r^d t$.

The \emph{polymer measure} $\DP_t^{\b,\om}$ is the Gibbsian probability measure on the space $\mathcal{C}([0,\infty)\times\mathbb{R})$ of continuous paths, defined by its density with respect to the Wiener measure:
\begin{equation}
\rmd \DP_t^{\b,\om} = \frac{1}{Z_t} \exp\left(\b \om(V_t)\right) \rmd \DP,
\end{equation}
where $\b\in\mathbb{R}^d$ is the inverse temperature parameter, and where $Z_t$ is the \emph{partition function} of the polymer\footnote{For every probability measure $\IP$ and variable $X$, we will use the convention that $\IP[X]$ denotes the expectation of $X$ under $\IP$.}:
\begin{equation} \label{eq:partition_function_def}
Z_t(\om,\b,r) = \DP\left[ \exp\left(\b \om(V_t)\right)\right].
\end{equation}
This model was first introduced by Comets and Yoshida \cite{CY05}. Under the polymer measure, the path is attracted by the Poisson points when $\b>0$, and repelled when $\b <0$.

The partition function $Z_t$ has mean $\exp(\lambda(\b) \nu r^d t)$, where
\begin{equation}
\lambda(\b) = e^\b - 1.
\end{equation}
The \emph{renormalized partition function}:
\begin{equation}
W_t(\om,\b,r) = e^{-\lambda \nu r^d t} Z_t(\om,\b,r),
\end{equation}
is a mean one, positive martingale such that $W_t\to W_\infty$. In \cite{CY05}, it is shown that there is a dichotomy:
\begin{equation}
\text{either} \quad W_\infty = 0 \  \text{ a.s.} \quad  \text{or} \quad W_\infty > 0 \  \text{ a.s.}
\end{equation}
The first case corresponds to the \emph{strong disorder regime} and the second to the \emph{weak disorder regime}. Furthermore, it is proved that there exist critical values $\b_{c}^{-}(\nu,r)\in [-\infty,0]$ and $\b_{c}^{+}(\nu,r)\in [0,\infty]$, such that weak disorder holds for  $\b_{c}^{-}<\b<\b_{c}^{+} $ and strong disorder holds for $\b\notin [\b_c^-,\b_c^+]$.

A classical quantity in statistical physics is the deterministic \emph{quenched free energy}:
\begin{equation}
p(\b,\nu,r) = \lim_{t\to\infty} \frac{1}{t} \ln Z_t(\om,\b,r),
\end{equation}
which can be compared to the quenched free energy, that is defined as the limit as $t\to \infty$ of $t^{-1} \ln \IP[Z_t]$ and equals $\lambda(\b) \nu r^d$. The fact that the two energies agree or not also separates two regions, with some critical $\b$ thresholds:
\begin{equation} \label{eq:very_strong}
\text{either} \quad p(\b,\nu,r) < \lambda(\b) \nu r^d \quad  \text{or} \quad p(\b,\nu,r) = \lambda(\b) \nu r^d.
\end{equation}
As strict inequality in \eqref{eq:very_strong} implies strong disorder, this regime is sometimes called \emph{very strong disorder}. When this holds, the difference of the two energies - also called the \emph{quenched Lyapunov exponent} - describes the exponential decay of $W_t$. Whether or not very strong disorder is equivalent to strong disorder is still an open question. 

Path localization properties of the polymer were studied by Comets and Yoshida in \cite{CYBMPO2}. For all time $t\geq 0$, they considered a \emph{favorite path} $Y^{(t)}_s$, depending only on the Poisson environment, such that for all $s\leq t$, \[\DP_t^{\b,\om}(B_s\in U(Y^{(t)}_s)) = \max_{x\in\mathbb{R}^d} \DP_t^{\b,\om}(B_s\in U(x)).\]
Denoting by $R_t^*(B)=\int_0^t \mathbf{1}
\{B_s \in U(Y^{(t)}_s)\}\rmd s$ the time fraction any path $B$ spends next to the favorite path, they showed that under strict inequality between the derivatives in $\b$ of the two free energies, the following localization property holds:
\begin{equation} \label{eq:loc_prop}
\liminf_{t\to\infty} \IP\DP_t^{\b,\om}[R_t^*] >0.
\end{equation}
When this is verified, the polymer spends $\DP_t^{\b,\om}$-a.s.\ignorespaces\  a positive fraction of time close to a particular path in the environment, which is in contrast with the Brownian motion behavior.

For a preliminary step, one can define a favorite endpoint for the polymer on $[0,t]$ in a similar way (see \cite[Remark 2.3.1]{CY05}). A complete understanding of the endpoint distribution has been recently achieved \cite{BatesChatterjee} in the discrete case.

Under weak disorder, the environment is supposed to have less influence over the polymer measure, so the polymer should behave similarly to Brownian motion ($\b = 0$) and paths should not be localized. In the case of the discrete polymer, there is a functional central limit theorem on the polymer path, which holds in the whole weak disorder region \cite{CY06}. For the Brownian polymer, a central limit theorem for the endpoint distribution was only proved in the smaller $L^2$ region \cite{CYkokyuroku}. In this region, a local limit theorem for the discrete and continuous polymers was proved by Vargas \cite{Va04}. One can observe that in general, continuous models have received less attention and results are still incomplete compared to the discrete models. In contrast with weak disorder, the strong disorder regime should be characterized by localized paths and super-diffusivity ($B_t \approx t^\xi$ as $t\to\infty$ with $\xi > 1/2$), but integrable models put aside, rigorous proofs of this facts seem still out of reach.

We end this section by mentioning two related models:
\begin{enumerate}[label=(\roman*)]
\item The branching Brownian motion in Poisson environment \cite{Shiozawa-clt,Shiozawa-loc}, for which the mean population size at a given time equals the partition function \eqref{eq:partition_function_def}.
\item The stochastic heat equation driven by a multiplicative Gaussian white noise with smoothing in space and rescaling, studied in \cite{Chiranjib-wkstrgdisorder}, in dimension $d\geq 3$.
\end{enumerate}

\subsection{KPZ universality for polymers and the intermediate disorder regime}
From now on, we focus on dimension $d=1$. In this case, the polymer is in the strong disorder phase as soon as $\b \neq 0$. It is expected that under the polymer measure,
\begin{equation} \label{eq:expected_critical_exponents}
\sup_{0\leq s \leq t} |B_s| \approx t^{2/3} \quad \text{and} \quad \ln Z_t - \IP[\ln Z_t] \approx t^{1/3} \  \text{ as} \  t\to\infty.
\end{equation}
Moreover, it is conjectured that:
\begin{conjecture} For all non-zero $\b,\nu$ and $r$, there exists a constant $\sigma(\beta,\nu,r)$ such that, as $t\to \infty$,
\begin{equation} \label{eq:fluctuationsExponents}
\frac{\ln Z_t-p(\b,\nu,r)t}{\sigma(\beta,\nu,r)t^{1/3}} \cvlaw F_{\mathrm{GOE}} \; ,
\end{equation}
where the $F_{\mathrm{GOE}}$ is the Tracy-Widom GOE distribution \textup{\cite{tracy1994level}}.
\end{conjecture}
These properties are characteristics of the KPZ universality class (cf. Section \ref{sec:kpz}). They are in sharp contrast to the weak disorder regime, where one knows to a large extent that $B_t \approx t^{1/2}$, and where the free energy $\ln Z_t$ has order one fluctuations around its mean \cite{CY05}.

For general models, only non-sharp bounds for the fluctuations in \eqref{eq:fluctuationsExponents} have been obtained yet \cite{CY05,Mejane,Petermann,Wuthrich-fluctuation,Wuthrich1}. For a specific, integrable and discrete model of polymer involving inverse log-gamma distributed weights and some boundary conditions, Seppäläinen \cite{Sepp12} was able to obtain sharp bounds with probabilistic methods, while the Tracy-Widom fluctuations were obtained later via Fredholm determinant identities \cite{BoCorFe14}.

A more accessible question is to look at the so-called \emph{intermediate disorder regime}, in the transition between $\b>0$ (strong disorder) and $\b=0$ (weak disorder). In the seminal paper \cite{AlbertsKhaninQuastelAP}, Alberts, Khanin and Quastel considered a time-space diffusive rescaling of the discrete directed polymer in i.i.d. \unskip \  environment, where they also rescaled the inverse temperature as $\b_n= \b n^{-1/4}$. They proved that under this scaling, the point-to-point and point-to-line partition functions of the polymer converge, in distribution, toward the partition functions of the \emph{continuum directed polymer}, which is a directed polymer of Brownian path and white noise environment. \medskip

We now outline the \textbf{main results in the paper} (Theorems \ref{th:cv_Wt_interReg} - \ref{th:cvlaw}). We prove that the intermediate disorder regime also appears as a scaling limit of the Brownian directed polymer in Poisson environment. Here, thanks to the Poissonian environment, the model has a more general scaling (compared to the general discrete model) that involves parameters $\b,\nu$ and $r$. In particular, at the price of tuning the other parameters, we show that the intermediate disorder regime can occur while keeping the temperature fixed. Similarly to \cite{AlbertsKhaninQuastelAP,PolynomialChaos17}, the result is obtained via chaos expansion of the partition functions. In our case,
$W_t$ admits an infinite Wiener-It\^o chaos expansion \cite{LastPenrose} which arises from the nice algebraic Poisson structure.

The paper is structured as follows:
In the rest of this introduction, we discuss the link between the P2P partition function of the polymer and the stochastic heat equation. We will also say a few words about the KPZ equation, the KPZ universality class and the intermediate disorder regime in the discrete setting. The main results are presented in Section \ref{sec:results}. Sections \ref{sec:wienerIto} and \ref{sec:wienerIntegrals} are devoted to introducing the chaos expansions in the Poisson and white noise environment. In Section \ref{sec:ustats}, we study the asymptotics of Poisson Wiener-It\^o integrals. Proofs of the results will finally be given through Section \ref{sec:proofs}.

\subsection{The KPZ equation and the stochastic heat equation}
The Kardar-Parisi-Zhang equation is the non-linear stochastic partial differential equation:\footnote{We will try to reserve capital letters ($T,X$,...) for the macroscopic scale (KPZ, SHE and continuum polymer) and lower case ($t,x$,...) for the microscopic scale (Poisson polymer and associated quantities).}
\begin{equation} \label{KPZ} \frac{\partial \mathcal{H}_\b}{\partial T}(T,X) = \frac{1}{2} \frac{\partial^2 \mathcal{H}_\b}{\partial X^2}(T,X) + \frac{1}{2} \left(\frac{\partial \mathcal{H}_\b}{\partial X}(T,X) \right)^2 + \beta \eta(T,X),
\end{equation}
where $\b\in\mathbb{R}$ and $\eta$ stands for the space-time Gaussian white noise (for a definition of this object, see section \ref{subseq:white_noise}).
The equation was first introduced in 1986 by Kardar, Parisi and Zhang \cite{KPZ}, in the study of scaling behaviors of
random interface growth.

Due to the non-linear term, it is difficult to give a proper definition of a solution of the KPZ equation. Bertini-Cancrini \cite{BertiniCancrini95} argued that $\mathcal{H}_\b$ could be defined by the so-called \emph{Hopf-Cole transformation}:
\begin{equation} \label{eq:Hopf-Cole}
\mathcal{H}_\b(T,X) = \ln \mathcal{Z}_\b(T,X),\end{equation}
where $\mathcal{Z}_\b$ is the solution of the stochastic heat equation (SHE):
\begin{equation} \label{SHE}
\frac{\partial \mathcal{Z}_\b}{\partial T}(T,X) = \frac{1}{2} \frac{\partial^2 \mathcal{Z}_\b}{\partial X^2}(T,X) + \beta \mathcal{Z}_\b(T,X)\eta(T,X).
\end{equation}
As a first-approach justification, one can check that the relation \eqref{eq:Hopf-Cole} defines a solution to \eqref{KPZ}, whenever $\mathcal{Z}$ is a solution of \eqref{SHE}, where the white noise $\eta$ is replaced with a smooth function.

Developing new tools to make sense of ill-posed stochastic PDEs, Hairer \cite{hairer} later constructed a method giving a direct notion of solution to the KPZ equation. Hairer further showed that the solution coincided with the solution defined by the Hopf-Cole transformation.

\subsection{Connections between the stochastic heat equation(s) and the directed polymers in random environments}
\subsubsection{The Poisson case}
Introduce the normalized \emph{point-to-point partition function}:
\begin{equation} \label{eq:renormalizedP2Pdef}
W(t,x;\om,\b,r) = \rho(t,x) \DP^{t,x}_{0,0}\left[ \exp \{ \b \om(V_t)-\lambda(\b) \nu r^d t\}\right],
\end{equation}
where $\DP^{t,x}_{0,0}$ is the Brownian bridge $(0,0)\to(t,x)$ and $\rho(t,x)=e^{-x^2/2t}/\sqrt{2\pi t}$ is the heat kernel/Brownian motion transition function.
In the next theorem, we state that the renormalized P2P partition function verifies a weak formulation of the following stochastic heat equation, with multiplicative Poisson noise:
\begin{equation}
\partial_t W(t,x) = \frac{1}{2} \Delta W(t,x) + \lambda W(t-,x) \bar{\om}(\rmd t \times U(x)).
\end{equation}
Let $\mathcal{D}(\mathbb R)$ denote the set of infinitely differentiable functions of compact support.
\begin{theorem}\label{th:she}
 For all $\varphi\in \mathcal{D}(\mathbb R)$ and $t\geq 0$, we have $\IP$-almost surely
\begin{multline} \label{eq:sheForPoisson}
\int_\mathbb{R} W(t,x)\varphi(x) \rmd x = \varphi(0) + \frac{1}{2}\int_0^t \rmd s \int_\mathbb{R} W(s,x) \Delta \varphi (x) \rmd x \\
+ \lambda \int_{\mathbb{R}} \rmd x \varphi(x) \int_{(0,t]\times \mathbb R} \bar{\om}(\rmd s, \rmd y) W(s-,x) \mathbf{1}_{|y-x|\leq r/2} \,.
\end{multline}
\end{theorem}
\noindent The proof of the theorem can be found in Section \ref{sec:proofShe}.
\subsubsection{The continuous case}
\label{continuumPolymerSHE} 
A special case of interest for the SHE, where $\mathcal{Z}_\b(T,X)$ can be interpreted as the point-to-point partition function of a directed polymer, positioned at $X=0$ at time $T=0$, is when
\begin{equation} \label{eq:ICofSHE}
\mathcal{Z}_\b(0,X) = \delta_0(X).
\end{equation}
In this case, $\mathcal{Z}_\b(T,X)$ can be expressed through the following shortcut (cf. Section \ref{subsec:constructionofP2P}):
\begin{equation} \label{eq:FeynmanKac0}
 \mathcal{Z}_\b(T,X) = \rho(T,X) \  \DP_{0,0}^{T,X} \left[ :\mathrm{exp}:\left(\beta \int_0^T \eta(u,B_u)\mathrm{d}u\right) \right],
\end{equation}
where $:\exp:$ denotes the Wick exponential.

Renormalization put aside, this equation is similar to the definition of the point-to-point partition of our polymer model \eqref{eq:renormalizedP2Pdef}. Alberts, Khanin and Quastel \cite{AlbertsKhaninQuastelCRP} were indeed able to construct a polymer measure, with P2P partition function given by $\mathcal{Z}_\b(T,X)$.
As both the environment (white noise) and the path (Brownian motion) of the polymer are continuous, it was named \emph{the continuum directed random polymer}.

Some special care has to be taken to construct this measure, as it can be shown that contrary to the Poissonian medium polymer, the continuum polymer measure is actually singular with respect to the Wiener measure (See \cite{AlbertsKhaninQuastelCRP}). A way to circumvent this issue is to define directly the transition functions of the polymer through the equation that they should satisfy.

At time horizon $T=1$, the transition functions of the polymer path $X$ are defined by
\begin{multline*}
\mathbb{P}_{1,\beta}^{\eta} (X_{t_1} \in \mathrm{d} x_1, \dots , X_{t_k} \in \mathrm{d}x_k) \\
= \frac{\mathcal{Z}_\b(t_k,x_k \, ; 1,*)}{\mathcal{Z}_\b(0,0 \, ; 1,*)} \prod_{j=0}^{k-1} \mathcal{Z}_\b(t_j,x_j \, ; t_{j+1},x_{j+1}) \, \mathrm{d}x_1\dots \mathrm{d}x_k \, ,
\end{multline*}
where $\mathcal{Z}_\b (S,Y;T,X)$ is the P2P partition function from $(S,Y)$ to $(T,X)$:
\[
\mathcal{Z}_\b (S,Y;T,X) = \mathcal{Z}_\b (T-S,X-Y)\circ \theta_{S,Y},
\] 
letting $\theta_{S,Y}$ be the shift by $(S,Y)$ in the white noise environment, and where 
\[\mathcal{Z}_\b(S,Y \, ; 1,*) = \int_\mathbb{R} \mathcal{Z}_\b(S,Y \, ; 1,X)\rmd X.\] The \emph{point-to-line partition function} of the continuum polymer at time horizon $T=1$ is given by
\begin{equation}
\mathcal{Z}_\b = \mathcal{Z}_\b(0,0 \, ; 1,*).
\end{equation}

\begin{remark}
One can check that the transition functions of the Brownian polymer in Poisson environment also satisfy the above equation.
\end{remark}

Similarly to the Poisson polymer, the P2P free energy $\mathcal{F}_\b(T,X)$ can be defined as
\begin{equation}
\mathcal{F}_\b(T,X) = \ln \frac{\mathcal{Z}_\b(T,X)}{\rho(T,X)},
\end{equation}
so that the free energy of the polymer and the solution of the KPZ equation are related by the equation:
\begin{equation}
\mathcal{F}_\b(T,X) = \mathcal{H}_\b(T,X) + X^2/2T + \ln \sqrt{2 \pi T}.
\end{equation}

\subsection{The KPZ universality class and the KPZ equation} \label{sec:kpz}

The KPZ equation belongs to a wide class of mathematical and physical models, called the KPZ universality class, which gathers models that share similar statistical behaviors under long time or scaling limits and particular scaling exponents (3-2-1 in time, space and fluctuations, as in \eqref{eq:expected_critical_exponents}). The reader may refer to \cite{corwin2016kardar} for a non-technical review on the KPZ universality class and the KPZ equation. Notable models, that are proven or conjectured to belong to this class, include paths in a random environment (directed polymers in random environment, first and last passage percolation), random growing interfaces (corner growth model), interacting particle systems (asymmetric simple exclusion process (ASEP)), stochastic PDEs and random matrices.

In 2011, the computation, by Amir, Corwin and Quastel \cite{ACQ11}, of the pointwise distribution of $\mathcal{H}_\b(T,X)$ was a breakthrough and provided the proof that the KPZ equation lied in the KPZ universality class. The result relied on two main results: the work of Bertini-Giacomin \cite{BertiniGiacomin97}, who were able to show that the solution of the KPZ equation appeared as the weak asymmetry limit of the ASEP height function; and on the papers of Tracy-Widom \cite{Tracy2008Erra,Tracy2008,Tracy2009}, in which the authors obtained formulas to express this height function distribution.

As a consequence of the ASEP weak asymmetry limit of Bertini-Giacomin, the KPZ equation can be seen as a crossover between the positive asymmetry ASEP (which belongs to the KPZ universality class) and the symmetric simple exclusion process (which belongs to the Edwards-Wilkinson universality class with 4-2-1 scaling). The \emph{weak KPZ universality conjecture} states that the KPZ equation is a universal object of the KPZ universality class. As a general idea, the KPZ equation should appear as a scaling limit at critical parameters for models that feature a phase transition between the Edward-Wilkinson class (4-2-1 scaling) and the KPZ class. In recent years, several techniques have been used to prove convergence of different models to the KPZ equation in weakly asymmetric regimes. The Hopf-Cole transformation is one of them and provides the opportunity to deal with the linear stochastic heat equation, instead of the non-linear KPZ equation; we will rely on it in this paper. This transformation may be applied to models which can be controlled after exponentiation, as for \cite{BertiniGiacomin97,
AlbertsKhaninQuastelAP,Dembo2016,
corwin2017,Labbe2017}. When one cannot apply the transformation, another tool is the pathwise approach introduced by Hairer \cite{hairer} and considered in \cite{hairer2017,Gubinelli2017,
HOSHINO2017}. Although very robust, the pathwise analysis requires quantitative estimates that may be hard to obtain. An alternative, relying on stationarity of the models, is the martingale approach developed in \cite{ GubinelliEnergysolutions,Gubinelli2013,
Goncalves2014} and applied to \cite{Diehl2017,Goncalves2014,goncalves2015,
Goncalves2017,Franco2016}.

\subsection{The intermediate disorder regime for the discrete polymer with i.i.d. \unskip \  weights}
The fact that the KPZ equation can emerge as a crossover regime appears in the result of Alberts, Khanin and Quastel \cite{AlbertsKhaninQuastelAP}, who showed that in dimension $d=1$, the rescaled logarithm of the point-to-point partition function of the discrete directed polymer (see below) converges to $\mathcal{H}_\b(T,X) = \ln \mathcal{Z}_\b(T,X)$. The P2P partition function of the polymer is defined as
\begin{equation}
Z_\b(n,x)=P(S_n=x)\times\DP[e^{\b\sum_{k=0}^n w(i,S_i)} | S_n = x],
\end{equation}
where $S$ is the simple symmetric random walk and $w(i,x)$ are i.i.d. \unskip \  random variables with finite exponential moments. They showed that, as $n\to\infty$,
\begin{equation} \label{eq:AlbertsKhaninQuastel}
\frac{\sqrt{n}}{2}Z_{\b n^{-1/4}} (nT,\sqrt{n}X)e^{-\mu(\b n^{-1/4})nT} \cvlaw \mathcal{Z}_\b (T,X),
\end{equation}
where $\mu(\b)=\ln \DP[e^{\b w(i,x)}]$, and where the limit in distribution is proven in terms of convergence of processes.

When $\b = 0$ (weak disorder), the polymer measure reduces to the simple symmetric random walk measure, so the model belongs to the Edwards-Wilkinson universality class. When $\b > 0$ (strong disorder), a long-standing conjecture is to prove that the discrete polymer model lies in the KPZ universality class. In particular, what is expected is the following, where the $1/3$ coefficient and the limiting Tracy-Widom distribution are characteristics of the KPZ universality class:
\begin{conjecture} \textup{\cite{AlbertsKhaninQuastelAP}} Suppose that the $w(i,x)$ are i.i.d. \unskip \  of finite fifth moment. Then, there exists some constants $c(\beta)$ et $\sigma(\beta)$ such that, as $n\to \infty$,
\begin{equation*}
\frac{\ln Z_{\b}(n,0)-c(\beta)n}{\sigma(\beta)n^{1/3}} \cvlaw F_{\mathrm{GUE}}.
\end{equation*}
\end{conjecture}

Because in the limit \eqref{eq:AlbertsKhaninQuastel}, $\b n^{-1/4} \to 0$, the KPZ equation - or equivalently the continuum directed random polymer - can be interpreted as a crossover regime between the weak disorder polymer regime and the strong disorder  regime, so that it was named in  \cite{AlbertsKhaninQuastelAP} the \emph{intermediate disorder regime}.
Moreover, the intermediate disorder regime features both characteristics of the strong disorder (a limiting universal law that does not depend on the law of the initial environment, a limiting polymer model that in the KPZ universality class) and the weak disorder (a diffusive scaling and a random local limit theorem for the endpoint density \textup{\cite{AlbertsKhaninQuastelAP}}).

\section{Main results} \label{sec:results}
We will show that a similar result to convergence \eqref{eq:AlbertsKhaninQuastel} holds in our model. In this perspective, we consider parameters $\beta_t\in\mathbb{R}$, $\nu_t > 0$ and $r_t > 0$ that depend on time $t$, and we fix a parameter $\b^*\in\mathbb{R}^*$. We introduce three asymptotic relations, when $t\to\infty$:
\begin{equation} \label{eq:paramRI}
\begin{gathered}
\text{(a)} \  \nu_t r_t^2 \lambda(\b_t)^2 \sim (\b^*)^{2}t^{-1/2}, \quad \text{(b)} \  \nu_t r_t^3 \lambda(\b_t)^3 \to 0,\\
\text{(c)} \  r_t/\sqrt{t} \to 0.
\end{gathered}
\end{equation}

\begin{remark}
Suppose the radius $r_t$ and the intensity $\nu_t$ to be constants. Then, the relations imply that $\beta_t$ scales like $t^{-1/4}$, as in  equation \eqref{eq:AlbertsKhaninQuastel}.
\end{remark}
\begin{remark}[Interpretation of \eqref{eq:paramRI} as disorder intensity]
In dimension $d\geq 3$, we see from the proof of Proposition 4.2.1 in \cite{CY05} that there exists a positive constant $c(d)$, such that the polymer lies in the $L^2$ region (which is a subregion of weak disorder) as soon as
$\nu r^{d+2} \lambda(\b)^2 < c(d)$.  Since conditions (a) and (c) imply that $\nu_t r_t^3 \lambda(\b_t)^2 \to 0$, the scaling for $d=1$ should here also be interpreted as a crossover between strong and weak disorder.
We say more about the scaling relations in Section \ref{subseq:ScalingRelations}.
\end{remark}
\begin{remark}
The relations can be compared to the regime of  \emph{complete localization} \cite{CYBMPO2}, corresponding to the extremal parameters regime ($r$ is fixed):
\begin{equation}
\nu \to \infty, \;|\beta|\leq \b_0, \quad {\rm such\ that} \quad \nu \beta^2 \to \infty.
\end{equation}
In the \emph{complete localization} regime, the polymer measure is highly concentrated around a favorite path, and the rest of the environment is neglected.

\end{remark}

\begin{theorem} \label{th:cv_Wt_interReg} Under conditions (a), (b), (c) and as $t\to\infty$:
\begin{equation}
W_t(\om^{\nu_t},\b_t,r_t) \cvlaw \mathcal{Z}_{\beta^*},
\end{equation}
where $\om^{\nu_t}$ is the Poisson point process with intensity measure $\nu_t \rmd s \rmd x$.
\end{theorem}
We will also show that the result extends to the renormalized point-to-point partition function:
\begin{equation}
W(s,y;t,x;\om,\b,r) = W(t-s,x-y;\om,\b,r)\circ \theta_{s,y},
\end{equation}
where $\theta_{s,y}$ denotes the shift of vector $(s,y)$ in the Poisson environment, i.e. $\theta_{s,y} (\sum_i \delta_{(s_i,y_i)}) = \sum_i \delta_{(s_i-s,y_i-y)}$.
\begin{theorem} \label{th:cv_P2PWt_interReg} Let $S,T\geq 0$ and $Y,X\in\mathbb{R}$. Under conditions (a), (b), (c) and as $t\to\infty$:
\begin{equation}
\sqrt{t} W\left(tS,\sqrt{t} \hspace{0.3mm}Y;t\hspace{0.3mm}T,\sqrt{t}\hspace{0.3mm}X;\om^{\nu_t},\b_t,r_t\right) \cvlaw \mathcal{Z}_{\beta^*}\left(S,Y;T,X\right).
\end{equation}
\end{theorem}
\begin{remark}
The $\sqrt{t}$ term appears here as a renormalization in the scaling of the heat kernel: $\sqrt{t}\rho\left(tT,\sqrt{t}X\right) = \rho(T,X)$.
\end{remark}

Let $\mathcal{D}'(\mathcal R)$ denote the space of distributions on $\mathbb{R}$, and $D\big([0,1],\mathcal{D}'(\mathcal R)\big)$ the space of c\`adl\`ag function with values in the space of distributions, equipped with the topology defined in \cite{mitoma1983tightness}.
We also define the rescaled and renormalized P2P partition function \eqref{eq:renormalizedP2Pdef}:
\begin{equation}
\mathcal{Y}_t\left(T,X\right) = \rho(T,X) W\left(tT,\sqrt{t}X;\om^{\nu_t},\b_t,r_t\right).
\end{equation}
The two variables function $\mathcal{Y}_t$ can be seen as an element of $D\left([0,1],\mathcal{D}'(\mathcal R)\right)$ through the mapping $\mathcal{Y}_t : T \mapsto \left(\varphi \mapsto \int \mathcal{Y}_t(T,X)\varphi(X)\rmd X\right)$. We have:
\begin{theorem} \label{th:cvlaw}
Suppose that $(\b_t)_{t\geq 0}$ is bounded by above. Then, as $t\to\infty$:
\begin{equation}
\mathcal{Y}_t \cvlaw \big(T\mapsto \mathcal{Z}_{\b^*}(T,\cdot)\big),
\end{equation}
where the convergence in distribution holds in $D\big([0,1],\mathcal{D}'(\mathcal R)\big)$.
\end{theorem}

%

\section{The Wiener-It\^o integrals with respect to Poisson process} \label{sec:WienerItoPoisson} \label{sec:wienerIto}
In this section, we expose the basic theory of multiple integration over Poisson processes. We rely on the reviews of G\"unter Last and Mathew Penrose \cite{LastPenrose,Last2016}.

Let us first introduce some notations that will prove useful throughout the paper.
For any $k\geq 1$, $s_1,\dots,s_k \in \mathbb{R}_+$ and $x_1,\dots,x_k \in \mathbb{R}$, write $\mathbf{s}=(s_1,\dots,s_k)$ and $\mathbf{x}=(x_1,\dots,x_k)$. Let 
\begin{equation}\Delta_k(u,t)=\{\mathbf{s}\in [u,t]^k \,|\, u<s_1<\dots<s_k\leq t\},
\end{equation}
be the $k$-dimensional simplex and $\Delta_k=\Delta_k(0,1)$. In addition, for any given function $g$ of $\mathbb{R}_+^k\times \mathbb{R}^k$, define the symmetrized version of $g$   
\begin{equation}
\mathrm{Sym} \, g  \, (\mathbf{s},\mathbf{x})=\frac{1}{k!} \sum_{\pi \in {\mathfrak S}_k}
 g(\pi \mathbf{s}, \pi \mathbf{x}),
\end{equation}
where ${\mathfrak S}_k$ denotes the set of permutation of $\{1,\dots,k\}$, and where any $\pi\in {\mathfrak S}_k$ acts on $\mathbb{R}^k$ by permutation of indices.
We say that a function is \emph{symmetric} whenever $\mathrm{Sym} \, g = g$. We see from the definition that the function $\mathrm{Sym} \, g $ is indeed symmetric.

\subsection{The factorial measures}
For all function $f$ on $\mathbb{R}_+\times \mathbb{R}$, define its sum over the Poisson process as $\om(f):=\iint f(s,x)\om(\rmd s \rmd x)$. Suppose for a moment that one wants to evaluate $\IP[\om(f)^2]$. A solution is to decompose the square of the sum as
\begin{equation}
\om(f)^2 = \sum_{(s,x)\in \om} f(s,x)^2 + \sum_{(s,x)\neq(s',x') \in\, \om} f(s,x) f(s',x'),
\end{equation}
in order to apply the multivariate Mecke equation \cite[Th. 4.4]{LastPenrose} 
to evaluate the expectation of the second term on the right-hand side. One can observe that this particular term is an integral over a certain point process measure, which depends on $\om$ and defines a measure on $\mathbb{R}_+^2\times \mathbb{R}^2$, called the second \emph{factorial measure} of $\om$. The first factorial measure is simply $\om$.

As an example, we saw that the simple functional $\om(f)^2$ can be written as a sum of two integrals over factorial measures. More generally, we will see in Theorem \ref{th:wienerItoChaos} that any square-integrable functionals of $\om$ can be expressed through an infinite sum of integrals over similar measures; in particular, it will be the case for the partition function $W_t$.

To give a proper definition of the factorial measures, let $t>0$ and let $\mathcal{B}_t$ denote the product Borel sets of $[0,t]\times\mathbb{R}$, $\,\mathcal{B}_t^{\otimes k}$ the product Borel sets of $[0,t]^k\times\mathbb{R}^k$.

\begin{definition}
For any positive integer $k$, define the $\mathbf k$\textbf{-th factorial measure} $\,\om_t^{(k)}$ to be the point process on $[0,t]^k\times\mathbb{R}^k$, such that, for any measurable set $\mathbf{A}\in \mathcal{B}_t^{\otimes k}$,
\begin{equation}
\om_t^{(k)}(\mathbf{A}) = \underset{{(s_1,x_1),\dots,(s_k,x_k) \in\, \om_t}}{{\sum}^{\neq}} \mathbf{1}_{((s_1,x_1),\dots,(s_k,x_k)) \in \mathbf{A}},
\end{equation}
where the sign $\neq$ indicates that the summation is made over pairwise different $(s_i,x_i)$.
Otherwise stated,
\begin{equation}
\om_t^{(k)} = \underset{{(s_1,x_1),\dots,(s_k,x_k) \in\, \om_t}}{{\sum}^{\neq}} \delta_{((s_1,x_1),\dots,(s_k,x_k))}\;.
\end{equation}
\end{definition}
\begin{remark}
If $A$ is a borel set of $\mathcal{B}_t$ and $A^k=A\times \ldots \times A$, then $\om_t^{(k)}(A^k)$ is the number of $k$-tuples of distinct points of $\om_t$ that belong to $A$, so that
\[\om_t^{(k)}(A^k) = \om_t(A)(\om_t(A)-1)\dots(\om_t(A)-k+1),\]
which is a reason why it is called a "factorial" measure.
For $A_1,\dots,A_k$ a collection of pairwise disjoint sets of $\mathcal{B}_t$, the situation is substantially different as we have
\begin{equation} \label{eq:fact_measure_fact}
\om_t^{(k)}(A_1\times\dots\times A_k) = \prod_{i=1}^k \om(A_i).
\end{equation}
\end{remark}
Since the sum is over all distinct $k$-tuples, symmetry plays an important role in factorial measures, and one should keep in mind that symmetric functions are the natural functions to integrate. As an example, the integral of a function $g$ is in fact the integral of its symmetrized function:
\begin{align*}
\om_t^{(k)}(\mathrm{Sym}\,g)
& = \frac{1}{k!} \sum_{\pi\in {\mathfrak S}_k} \underset{{(s_i,x_i) \in\, \om_t}}{{\sum}^{\neq}}  g(\pi \mathbf{s}, \pi \mathbf{x})\\
& =\frac{1}{k!} \sum_{\pi\in {\mathfrak S}_k} \underset{{(\tilde{s}_i,\tilde{x}_i) \in\, \om_t}}{{\sum}^{\neq}}  g(\mathbf{\tilde{s}}, \mathbf{\tilde{x}})\\
& = \om_t^{(k)}(g).
\end{align*}

\subsection{Multiple stochastic integral over a Poissonian medium}
Now that we have defined the factorial measures of $\om_t$, we seek to do the same for the compensated measure $\bar{\om_t}$ - also called the \emph{Wiener-It\^o integral}. In particular, we still want to avoid points belonging to the diagonal in the integration process.

Let $A_1, \dots A_k$ be a collection of pairwise disjoint, finite sets of $\mathcal{B}_t$. Then, observe that
\begin{align*}
\prod_{i=1}^k \bar{\om_t}(A_i)&=\prod_{i=1}^k \big({\om_t}(A_i)-\nu|A_i|\big)\\
&=\sum_{J\subset [k]} \prod_{i\in J} \left(\int \mathbf{1}_{A_i}(s,x) \, \om_t(\rmd s \rmd x) \right) \prod_{i\in J^c} \left( -\int \mathbf{1}_{A_i}(s,x) \, \nu \rmd s \rmd x \right),
\end{align*}
where $[k] = \{1,2,\dots,k\}$ will not be confused with the integer part. Using the fact that the $A_i$ are disjoint, by \eqref{eq:fact_measure_fact} the above product over $J$ can be written as an integral with respect to the measure $\om^{(|J|)}_t$:
\begin{equation} \label{eq:formule_produit}
\begin{split}
&\prod_{i=1}^k \bar{\om_t}(A_i)\\
&=\sum_{J\subset [k]} (-1)^{k-|J|} \int_{[0,t]^k\times \mathbb{R}^k} \left(\prod_{i=1}^k\mathbf{1}_{A_i}\right) \om_t^{(|J|)}(\rmd \mathbf{s}_J,\rmd \mathbf{x}_J) \, \nu^{k-|J|} \,\rmd \mathbf{s}_{J^c}\,\rmd \mathbf{x}_{J^c},
\end{split}
\end{equation}
where the notations $\rmd \mathbf{s}_J$ and $\rmd \mathbf{x}_J$ mean that the integration is done with respect to the variables $(s_i)_{i\in J}$ and $(x_i)_{i\in J}$.
This leads to the following definition of the \emph{multiple Wiener-It\^o integral}:

\begin{definition}
For $k\geq1$ and $g\in L^1([0,t]^k\times \mathbb{R}^k)$, denote
the \textbf{multiple Wiener-It\^o integral} of $g$ as
\begin{equation} \label{eq:def_omk}
\begin{split}
&\bar{\om}^{(k)}_t(g)\\
&:=\sum_{J\subset [k]} (-1)^{k-|J|} \int_{[0,t]^k\times \mathbb{R}^k} 
g( \mathbf{s}, \mathbf{x})\, \om_t^{(|J|)}(\rmd \mathbf{s}_J,\rmd \mathbf{x}_J) \, \nu^{k-|J|} \,\rmd \mathbf{s}_{J^c} \, \rmd \mathbf{x}_{J^c}.
\end{split}
\end{equation}
When $k=0$, define $\bar{\om}^{(0)}_t$ to be the identity on $\mathbb{R}$.
\end{definition}
%

The two following results can be found in \cite{LastPenrose,Last2016}:
\begin{proposition} 
For $k\geq 1$, the map $\bar{\om_t}^{(k)}$ can be extended to a map \[\begin{array}{ccccc}
\bar{\om}_t^{(k)} & : & L^2([0,t]^k\times \mathbb{R}^k) & \to & L^2(\Omega,\mathcal{G},\IP)\\
& & g &\mapsto & \bar{\om}_t^{(k)}(g),
\end{array} \] which coincides with the above definition of $\,\bar{\om}_t^{(k)}$ on the functions of $L^1 \cap L^2([0,t]^k\times \mathbb{R}^k)$.
\end{proposition}
\begin{property}
\begin{enumerate}[label=(\roman*)]
\item For $k\geq 1$ and $g\in L^2([0,t]^k\times \mathbb{R}^k)$, then $\IP\big[\bar{\om}_t^{(k)} (g)\big] = 0$.
\item For every square-integrable function $g$,
\begin{equation}\label{eq:sym_om} \bar{\om_t}^{(k)}(\mathrm{Sym} \, g) = \bar{\om_t}^{(k)}(g).
\end{equation}
\item For any $k\geq 1$ and $l\geq 1$, $g\in L^2([0,t]^k\times \mathbb{R}^k)$ and $h\in L^2([0,t]^l\times \mathbb{R}^l)$, the following covariance structure holds:
\begin{equation} \label{eq:covariance_structure}
\IP\left[\bar{\om}_t^{(k)}(g) \, \bar{\om}_t^{(l)}(h) \right] = \delta_{k,l} \  k! \, \nu^k <\mathrm{Sym} \, g ,\mathrm{Sym} \, h >_{L^2([0,t]^k\times \mathbb{R}^k)}.
\end{equation}
\item The map $\bar{\om}_t^{(k)}$ is linear, in the sense that for all square-integrable $f,g$ and reals $\lambda,\mu$,
\begin{equation*}
\IP \text{-a.s.} \quad \bar{\om}_t^{(k)}(\lambda f + \mu g) = \lambda \, \bar{\om}_t^{(k)}(f) + \mu \,\bar{\om}_t^{(k)}(g).
\end{equation*}
\end{enumerate}
\end{property} 
{\bf Multiple Wiener-It\^o integral of a function on the simplex:} When $g$ is a function of $L^2(\Delta_k(0,t)\times\mathbb{R}^k)$, denote by $\widehat{g}$ its extension set to zero outside of the simplex, and define
\begin{equation} \label{eq:XWienerIto}
\bar{\om}_t^{(k)}({g}) := \bar{\om}_t^{(k)}(\widehat{g}).
\end{equation}
Observing that $[0,t]^k\times \mathbb{R}^k$ is made of $k!$ copies of $\Delta_k(0,t)\times \mathbb{R}^k$, one gets that
\begin{equation} \label{eq:iso_simpl}
\IP\left[{\bar{\om}_t^{(k)}({g})}^2 \right]= \nu^k \Vert g\Vert^2_{L^2(\Delta_k(0,t)\times\mathbb{R}^k)}.
\end{equation}
\begin{remark}
This tells us that $\bar{\om}_t^{(k)}$ is an isometry on $L^2(\Delta_k(0,t)\times\mathbb{R}^k)$ with measure $\nu^k\rmd \mathbf{s} \rmd \mathbf{x}$. This is one of the reasons why we will mainly consider functions of the simplex.
\end{remark}

\subsection{A Wiener-It\^o Chaos Expansion of the renormalized partition function} \label{subsection_WIchaos_expansion}
The next theorem, proven in \cite[\S 18.4]{LastPenrose}, states that every square-integrable function $F(\om)$ admits a Wiener-It\^o chaos expansion, i.e., it can be written as an infinite sum of orthogonal multiple Wiener-It\^o integrals. In order to be able to state the theorem, we first need to introduce a new operator.

For any function $F(\om_t)$ of the point-process up to time $t$, define the derivative operator $D_{(s,x)}(F) = F(\om_t + \delta_{(s,x)}) - F(\om_t)$, for all $(s,x)\in[0,t]\times\mathbb{R}$. Now, for any given $(s_i,x_i)\in[0,t]\times\mathbb{R}$, $i\leq k$, define the iterated operator:
\begin{eqnarray}
D_{(s_1,x_1),\dots,(s_k,x_k)}&=& D_{(s_1,x_1)}\circ D_{(s_2,x_2)} \circ\dots\circ D_{(s_k,x_k)} \, , \end{eqnarray}
and the function $T_k F : [0,t]^k \times \mathbb{R}^k \to \mathbb{R}$, by letting
\begin{equation}
T_k F (s_1,\dots,s_k,x_1,\dots,x_k) = \IP \big[D_{(s_1,x_1),\dots,(s_k,x_k)} F(\om_t) \big].
\end{equation}
We also set $T_0 F = \IP[F(\om_t)]$. By induction, we see that
$$
D_{(s_1,x_1),\dots,(s_k,x_k)} F (\omega_t) = \sum_{I \subset [k]} (-1)^{k-|I|} F\left( \omega_t + \sum_{i \in I} \delta_{(s_i,x_i)}\right) \,. 
$$
\begin{theorem} \label{th:wienerItoChaos}
Let $F(\om_t)$ be any measurable function of the point process up to time $t$, verifying $\IP\big[F(\om_t)^2\big]<\infty$. Then, for all $k$, $T_k F$ is a symmetric, square-integrable function, and we have $\IP$-almost surely:
\begin{equation} \label{eq:chaos_dec}
F(\om_t)=\sum_{k=0}^\infty \frac{1}{k!} \bar{\om}^{(k)}_t \big(T_k F\big).
\end{equation} 
The orthogonal series  \eqref{eq:chaos_dec} converges in $L^2(\Omega, \mathcal{G}, \P)$, and it is called the \textbf{Wiener-It\^o chaos expansion} of $F(\om_t)$.
\end{theorem}

\begin{remark} Note that the terms in the sum are pairwise orthogonal. 
As a possible interpretation of the theorem, one can view $T_k$ as a $k$-th derivative, and the Wiener-It\^o expansion as a Taylor expansion. 
\end{remark}

We can now apply the theorem to renormalized partition function $W_t$, which is square-integrable. For this purpose, define for all path $B$ and all $\delta>0$,
\begin{equation}
\chi^\delta_{s,x}(B) = \mathbf{1} \{|x-B_s|\leq\delta/2\}.
\end{equation}
Note that when $\delta = r$, we have $\chi^r_{s,x}(B) = \mathbf{1} \{x\in U(B_s)\}$.

\begin{proposition} \label{prop:Poisson_chaos_decomposition}
The renormalized partition function admits the following Wiener-It\^o chaos expansion:
\begin{equation}
W_t = \sum_{k=0}^\infty \frac{1}{k!} \, \bar{\om}^{(k)}_t \left(T_k \, W_t\right),
\end{equation} where, for all $\mathbf{s}\in [0,t]^k$, $\mathbf{x}\in\mathbb{R}^k$ and $k\geq 0$:
\begin{equation} \label{eq:TkWtExpression}
T_k \, W_t \,(\mathbf{s},\mathbf{x}) = \lambda(\b)^k\, \DP\left[\prod_{i=1}^k \chi^r_{s_i,x_i}(B) \right],
\end{equation}
with the convention that an empty product equals $1$.
\end{proposition}

\begin{proof}
Fix a path $B$ and observe that  
\[\om\left(V_t(B)\right) = \int_{[0,t]\times\mathbb{R}} \chi^r_{s,y}(B)\, \om(\rmd s \rmd y).\] For all $(s,x)\in [0,t]\times \mathbb{R}$, we get that $D_{s,x}\, e^{\b \om(V_t)} = \lambda(\b)\,\chi^r_{s,x} \,e^{\b \om(V_t)} $. Hence,

\[T_k \, W_t \,(\mathbf{s},\mathbf{x}) = \IP \DP\left[e^{\b \om(V_t)} \prod_{i=1}^k \lambda \, \chi^r_{s_i,x_i}  \right] e^{-t \nu \lambda r^d}= \lambda ^k \, \DP\left[\prod_{i=1}^k \chi^r_{s_i,x_i} \right].\]
\end{proof}

\section{The Wiener integrals} \label{sec:wienerIntegrals}
We repeat the construction of Section \ref{sec:WienerItoPoisson}
for white noise instead of Poisson noise. 
In this section, we consider another probability space $(\Lambda,\mathcal{F}_\eta,\mathbb{Q})$.
\subsection{Stochastic integral over the white noise}
\label{subseq:white_noise}
\begin{definition}
A time-space \textbf{Gaussian white noise} environment $\eta$, is a random measure on $[0,1]\times\mathbb{R}$, which satisfies the two following properties:
\begin{enumerate}[label=(\roman*)]
\item For all measurable sets $A_1,\dots,A_k$ of $[0,1]\times\mathbb{R}$, $\big(\eta(A_1),\dots,\eta(A_k)\big)$ is a centered Gaussian vector.
\item For all measurable sets $A,B$ of $\,[0,1]\times\mathbb{R}$,
\begin{equation}
\mathbb Q[\eta(A)\eta(B)] = |A\cap B|.
\end{equation}
\end{enumerate}
\end{definition}
In the following, we suppose that a white noise process is defined on the probability space $(\Lambda,\mathcal{F}_\eta,\mathbb Q)$. It is then possible to construct a stochastic integral over the white noise measure, which has the following properties:

\begin{proposition}
There exists an isometry $I_1 : L^2\big([0,1]\times\mathbb{R}\big) \mapsto L^2(\Lambda,\mathcal{F}_\eta,\mathbb Q)$ verifying that:
\begin{enumerate}[label=(\roman*)]
\item For all measurable set $A$ of $[0,1]\times\mathbb{R}$, we have $I_1(A) = \eta(A)$.
\item For all $g\in L^2$, the variable $I_1(g)$ is a centered Gaussian variable of variance $\Vert g \Vert^2_{L^2([0,1]\times \mathbb{R})}$.
\end{enumerate}
\end{proposition}
We call $I_1(g)$ the \emph{stochastic integral} of $g$ \emph{over the white noise}. Note that the integral will sometimes also be written as
\begin{equation}
I_1(g) = \int_{[0,1]}\int_{\mathbb{R}} g(s,x) \, \eta(\rmd s,\rmd x).
\end{equation}

\subsection{Multiple stochastic integral}
It is again possible to extend $I_1$ to a multiple stochastic integral. One can find the details of such a procedure in Janson's book \cite[Chapter 7]{janson_1997}. This integral has very similar properties to the Wiener-It\^o integral:
\begin{theorem}  For all $k>0$, there exists a map $I_k:L^2([0,1]^k \times \mathbb{R}^k)\mapsto L^2(\Lambda,\mathcal{F}_\eta,\mathbb Q)$, which has the following properties:
\begin{enumerate}[label = (\roman*)]
\item If $g$ is any square-integrable function, then
$I_k(\mathrm{Sym} \, g) = I_k(g).$
\item For all $g\in L^2([0,1]^k \times \mathbb{R}^k)$ and $h\in L^2([0,1]^j \times \mathbb{R}^j)$, the variable $I_k(g)$ is centered and
\begin{equation}
 \mathbb Q \left[{I}_k(g) {I}_j(h)\right] =\ \delta_{k,j} \, k! \, <\mathrm{Sym} \, g,\mathrm{Sym} \, h>_{L^2([0,1]^k \times \mathbb{R}^k)}.
\end{equation}
\item For all orthogonal family $(g_1,\dots,g_k)$ of functions in $L^2([0,1] \times \mathbb{R})$, we have 
\begin{equation} \label{eq:IkOrthog} {I}_k\left(\prod_{j=1}^k g_j \right) =  \prod_{j=1}^k I_1(g_j).
\end{equation}
\item The map $I_k$ is linear.
\end{enumerate}
\end{theorem}
\begin{remark} \label{rk:simplex}
Similar to \eqref{eq:XWienerIto}, we define multiple Wiener integral of a function defined on the simplex.
 If $g$ is a function of $L^2(\Delta_k\times \mathbb{R})$, and if $\widehat{g}$ is the extension of $g$ set to zero outside of the simplex, we define $I_k(g) := I_k(\widehat g)$ and have
\begin{equation}
\mathbb Q[I_k(g)^2] = \Vert g \Vert^2_{L^2(\Delta_k\times\mathbb{R}^k)}.
\end{equation}
\end{remark}
\begin{remark}
We will sometimes use the notation
\begin{equation} I_k(g) = \int_{[0,1]^k} \int_{\mathbb{R}^k} g(\mathbf{t},\mathbf{x}) \eta^{\otimes k}(\mathrm{d}\mathbf{t},\mathrm{d}\mathbf{x}),
\end{equation}
where $[0,1]^k$ can be replaced with $\Delta_k$ when dealing with functions of the simplex. See \cite{janson_1997} for a justification of the tensor product notation.
\end{remark}

\subsection{Wiener chaos decomposition}
\begin{definition} 
For any family $G=(g^k)_{k\geq 0}$ such that for all $k\geq 0$, $g^k \in L^2(\Delta_k\times\mathbb{R}^k)$, and that
\[\Vert G \Vert^2_2:=\sum_{k=0}^{\infty} \Vert g^k \Vert^2_{L^2(\Delta_k\times\mathbb{R}^k)} < \infty,\]
we say that $G$ is an element of the \textbf{Fock space} $\bigoplus_{k=0}^\infty L^2(\Delta_k\times\mathbb{R}^k)$, which is a normed vector space with norm $\Vert G\Vert_2$, also called \textbf{Fock norm}.
\end{definition}

The next proposition is a consequence of Remark \ref{rk:simplex}:
\begin{proposition} The linear map
\begin{equation}
\begin{array}{ccccc}
I & : & \bigoplus_{k=0}^\infty L^2(\Delta_k\times\mathbb{R}^k) & \to & L^2(\Omega,\mathcal{F}_\eta,\mathbb Q)\\
& & G=(g^k)_{k\geq 0} &\mapsto & \sum_{k=0}^{\infty} I_k(g^k)=:I(G),
\end{array}
\end{equation}
is an isometry. Note that the sum is well defined as an $L^2-$limit.
\end{proposition}
\begin{remark}
If we were dealing with functions of $[0,1]^k\times\mathbb{R}^k$, we would have defined $I(G)$ as $\sum_{k=0}^\infty \frac{1}{k!} I_k(g^k)$ and the Fock norm as $\sum_{k=0}^\infty \frac{1}{k!} \Vert g^k \Vert^2_2$.
\end{remark}

We now consider a key example. Let $k>0$, and introduce the $k$-th dimensional Brownian transition function, for $(\mathbf{s},\mathbf{x}) \in \Delta_k\times\mathbb{R}^k$:
\begin{equation}
\begin{split}
\rho^k(\mathbf{s},\mathbf{x}) &= \DP\left(B_{s_1} \in \rmd x_1,\dots, B_{s_k} \in \rmd x_k\right)\\
& = \rho(s_1,x_1) \left( \prod_{j=1}^{k-1} \rho( s_{j+1}-s_j,x_{j+1}-x_j) \right),
\end{split}
\end{equation}
with the convention that $\rho^0 = 1$.

\begin{proposition} \label{prop_rhoIsInFock} The family $R(\beta) = (\beta^k \rho^k)_{k\geq 0}$ is in $\bigoplus_{k=0}^\infty L^2(\Delta_k\times\mathbb{R}^k)$ for all $\b \in\mathbb R$. In particular, the variable
\begin{equation} \label{eq:defOfP2LZBeta}
\mathcal{Z}_\beta := I\big(R(\beta)\big),
\end{equation} 
is well defined and square-integrable.
\end{proposition}
\begin{remark}
This quantity is in fact the continuum polymer partition function (cf. Section \ref{subsec:constructionofP2P}).
\end{remark}
\begin{proof} We will rely on the observation that
\[\rho(s,x)^2 = \frac{1}{2\sqrt{\pi s}} \rho(s/2,x).\]
Expressing $\rho^k$ in terms of product of $\rho$ function, we have, with the convention that $s_0 = x_0 = 0$,
\[ \rho^k(\mathbf{s},\mathbf{x})^2 = \rho^k(\mathbf{s}/2,\mathbf{x}) \prod_{j=1}^k \frac{1}{2\sqrt{\pi(s_j-s_{j-1})}},\]
so,
\[\int_{\Delta_k} \int_{\mathbb{R}^k} \rho^k(\mathbf{s},\mathbf{x})^2 \mathrm{d}\mathbf{s} \mathrm{d}\mathbf{x}=2^{-k}\pi^{-k/2}\int_{\Delta_k} \prod_{j=1}^k \frac{1}{\sqrt{s_j-s_{j-1}}}\rmd \mathbf{s}.\]
The last integral is the normalizing constant of the order $k+1$ Dirichlet distribution, taken with parameter $\mathbf{\alpha}=(\frac{1}{2},\dots,\frac{1}{2},1)$. As this constant is known to be the multivariate Beta function $\prod_{i=1}^{k+1} \Gamma(\alpha_i) \Gamma\big(\sum_{i=1}^{k+1} \alpha_i\big)^{-1}$, we obtain, using the identity $\Gamma(1/2) = \sqrt{\pi}$, that
\begin{equation*}
\Vert R(\beta)\Vert_2^2 = \sum_{k=0}^\infty \frac{(\beta^k)^2}{2^k \Gamma(k/2+1)}<\infty.
\end{equation*}
\end{proof}

\subsection{Construction of the P2P and P2L functions of the continuum polymer}
\label{subsec:constructionofP2P}
It is said that $\mathcal Z$ is a \emph{mild} solution to the stochastic heat equation \eqref{SHE} if, for all fixed $0 \leq S <T\leq 1$,
\begin{equation} \label{eq:mildSHE}
\begin{split}
\mathcal{Z}(T,X) = & \int_{\mathbb{R}} \rho(T-S,X-Y)\mathcal{Z}(S,Y)\mathrm{d}Y \\
&+\beta \int_S^T \int_{\mathbb{R}} \rho(T-U,X-Y) \mathcal{Z}(U,Y) \eta(U,Y) \mathrm{d}U \mathrm{d}Y,
\end{split}
\end{equation}
and if for all $T\geq 0$, $\mathcal{Z}(T,X)$ is measurable with respect to the white noise on $[0,T]\times\mathbb{R}$.

\begin{remark}
As a motivation to look at this form of the equation, one can check that if $\mathcal{Z}(T,X)$ satisfies \eqref{eq:mildSHE} with a smooth deterministic function $\eta(U,Y)$, then $\mathcal{Z}(T,X)$ is a solution to the SHE \eqref{SHE} with smooth noise.
\end{remark}

\begin{remark}
Under some integrability condition, it can be shown that there is a unique mild solution - up to indistinguishability - to the SHE with Dirac initial condition \textup{\cite{BertiniCancrini95}}. This solution is continuous in time and space for $(T,X)\in (0,1]\times \mathbb{R}$, and it is continuous in $T=0$ in the space of distributions. Furthermore, $\mathcal{Z}_\b(T,X)$ can be shown to be positive for all $T>0$ \textup{\cite{Moreno14,Mueller}}.
\end{remark}

Using the initial condition $\mathcal{Z}(0,X) = \delta_X$, we get by iterating equation \eqref{eq:mildSHE} for $S=0$, that
\begin{align*}
&\mathcal{Z}(T,X) = \rho(T,X) + \beta \int_0^T \int_\mathbb{R} \rho(T-U,X-Y) \rho(U,Y) \eta(U,Y) \, \mathrm{d}U \mathrm{d}Y \\
&+ \beta^2 \iint_{0<R<U\leq T} \iint_{\mathbb{R}^2} \rho(T-U,X-Y) \rho(U-R,Y-Z) \mathcal{Z}(R,Z) \\
& \qquad \qquad \qquad \qquad \qquad \times \eta(U,Y) \eta(R,Z) \, \mathrm{d}U \mathrm{d}Y  \mathrm{d}R\, \mathrm{d}Z.
\end{align*}
By repeating this iteration procedure, the following expansion arises:
\begin{equation} \label{eq:ChaosExpantionSHE}
\mathcal{Z}(T,X) = \sum_{k=0}^\infty \b^k I_k\big( \rho^k(\cdot ; 0,0;T,X) \big),
\end{equation}
where we have used the notation, for $\mathbf{s}\in \Delta_k(s,t)$ and $x\in\mathbb{R}^d$,
\begin{equation}
\begin{split}& \rho^k(\mathbf{s},\mathbf{x} \, ; \, s,y \, ; t,x)\\
& = \rho(s_1-s,x_1-y) \left( \prod_{j=1}^{k-1} \rho(s_{j+1}-s_j,x_{j+1}-x_j) \right) \rho(t-s_k,x-x_k),
\end{split}
\end{equation}
with the convention that $\rho^0(\cdot,\cdot\,;s,y;t,x) = \rho(t-s,x-y)$.

The ratio $ \frac{\rho^k(\mathbf{s},\mathbf{x} \, ; \, s,y \, ; t,x)}{\rho( t-s,x-y)} $ is the $k$-steps transition function of a Brownian bridge,  starting from $(s,y)$ and ending at $(t,x)$. From this observation, it is possible to introduce an alternative expression of $\mathcal{Z}(T,X)$, via a Feynman-Kac formula:
\begin{equation} \label{eq:FeynmanKac}
 \mathcal{Z}(T,X) = \rho(T,X) \  \DP_{0,0}^{T,X} \left[ :\mathrm{exp}:\left(\beta \int_0^T \eta(u,B_u)\mathrm{d}u\right) \right],
\end{equation}
with $\DP_{s,y}^{t,x}$ the Brownian bridge $(s,y) \to (t,x)$. The Wick exponential $:\mathrm{exp}:$ of a Gaussian random variable $\xi$ is defined by 
\[:\mathrm{exp}(\xi): = \sum_{k=0}^\infty \frac{1}{k!} :\xi^k:
 \]
where the $:\xi^k:$ notation stands for the Wick power of a random variable (cf.\cite{janson_1997}).
The integral $\int_0^T \eta(u,B_u)\mathrm{d}u$, on the other hand, is not well defined, and to understand how to go from \eqref{eq:FeynmanKac} to  \eqref{eq:ChaosExpantionSHE}, one should use the following shortcut:
\[ \DP_{0,0}^{T,X} \left[ :\left( \b \int_0^T \eta(u,B_u)\mathrm{d}u \right)^k : \right]  = \beta^k k! \int_{\Delta_k} \int_{\mathbb{R}^k} \frac{\rho^k(\mathbf{t},\mathbf{x}; \, 0,0 \, ; T,X)}{\rho(T,X)} \eta^{\otimes k}(\mathrm{d}\mathbf{t} \mathrm{d}\mathbf{x})\]

We will from now on suppose that $\mathcal{Z}_\b(T,X)$ is defined through equation \eqref{eq:ChaosExpantionSHE}. Again, it is possible to check  that the $\rho^k(\cdot \, ;\, S,Y;T,X)$ have finite Fock space norm \cite{AlbertsKhaninQuastelCRP}. One can further define the shifted P2P functions:
\begin{equation} \label{eq:P2Pdefinition} \mathcal{Z}_\b (S,Y;T,X) = \sum_{k=0}^\infty \b^k I_k\big( \rho^k(\cdot \, ;\, S,Y;T,X) \big).
\end{equation}
Moreover, integrating over $X$ equation \eqref{eq:ChaosExpantionSHE}, we recover the previous definition of the P2L partition function \eqref{eq:defOfP2LZBeta}:
\begin{equation*}
\mathcal{Z}_\b = \mathcal{Z}_\b(0,0;1,*) = \sum_{k=0}^\infty \b^k I_k( \rho^k)=I(R(\beta)).
\end{equation*}
We also get that for any test function $\varphi \in \mathcal{C}^\infty_c$:
\begin{equation} \label{eq:sheContreTestfunction}
\int_\mathbb{R} \mathcal{Z}_\b(T,X) \varphi(X) \rmd X = \sum_{k=0}^\infty \b^k I_k\left( \int \rho^k(\cdot \,; 0,0;T,X) \varphi(X) \rmd X \right).
\end{equation}

\section{Asymptotic study of Wiener-It\^o integrals} \label{sec:ustats}

\subsection{The scaling relations} \label{subseq:ScalingRelations}
From now on, we will suppose that, as $t\to\infty$,
\begin{equation}
\begin{gathered}
\text{(a)} \  \nu_t r_t^2 \lambda(\b_t)^2 \sim (\b^*)^{2}t^{-1/2}, \quad \text{(b)} \  \nu_t r_t^3 \lambda(\b_t)^3 \to 0,\\
\text{(c)} \  r_t/\sqrt{t} \to 0.
\end{gathered}
\end{equation}

There are two main reasons why we chose these relations. First, as $t\to\infty$, conditions (a) and (b) assure that under a scaling of $t$ in time and $\sqrt{t}$ in space, the Poisson environment over time $[0,t]\times\mathbb{R}$ becomes a white noise environment on $[0,1]\times\mathbb{R}$. This fact is properly stated in Theorem \ref{thm_ustat}, and the $k=1$ case of the proof gives good insights about how the parameters relate to one another. The addition of condition (c) ensures that the properly rescaled and normalized $T_k W_t$ functions converge to the Brownian transition functions.
\subsection{Gaussian limits of Wiener-It\^o integrals}
We are interested in the limit of rescaled Wiener-It\^o integrals, and more generally at sums of these integrals. We show that we can adapt the techniques developed for the study of U-statistics made out of an i.i.d. \unskip \  sequence of random variables, in chapter 11 of \cite{janson_1997}.

For any function $g$ defined on $[0,1]^k\times\mathbb{R}^k$ (resp. $\Delta_k\times\mathbb{R}^k$), denote by $\tilde{g}_t$ the rescaled function of $g$, defined on $[0,t]^k \times \mathbb{R}^k$ (resp. $\Delta_k(0,t)\times\mathbb{R}^k$) and such that 
\begin{equation}
\tilde{g}_t(\mathbf{s},\mathbf{x}) = g\big(\mathbf{s} / t,\mathbf{x} / \sqrt{t}\big),
\end{equation}
and let $\gamma_t$ be proportional to the vanishing parameter appearing in (b):
\begin{equation}
\gamma_t := (\b^*)^{-3} \nu_t r_t^3 \lambda(\b_t)^3 \to 0.
\end{equation}
\begin{theorem} \label{thm_ustat}
Let $g\in L^2([0,1]^k\times \mathbb{R}^k)$ for $k\geq 0$. The following convergence holds, as $t\to \infty$,
\begin{equation} 
\gamma_t^{k} \, \bar{\om}_t^{(k)}(\tilde{g}_t) \cvlaw \int_{[0,1]^k} \int_{\mathbb{R}^k} g(\mathbf{t},\mathbf{x}) \eta^{\otimes k}(\mathrm{d}\mathbf{t},\mathrm{d}\mathbf{x}).
\end{equation}

The convergence can be extended for any finite collection  $k_1,\dots ,k_m\in \mathbb{N}$ and $g^1, \dots ,g^m$, which satisfy $g^i\in L^2([0,1]^{k_i}\times \mathbb{R}^{k_i})$,
\begin{equation} \label{cvloi_th_2}
\big(\gamma_t^{k_1} \, \bar{\om}_t^{(k_1)}(\tilde{g}^{\,1}_t),  \dots  ,\gamma_t^{k_m}\,\bar{\om}_t^{(k_m)}(\tilde{g}^{\,m}_t) \big)
\cvlaw \big(I_{k_1}(g^1),\dots ,I_{k_m}(g^m)\big).
\end{equation}
\end{theorem}

\begin{corollary} \label{coro_ustat}
Let $G=(g^k)_{k\geq 0}$ belong to the Fock space $\bigoplus_{k\geq 0} L^2(\Delta_k\times \mathbb{R}^k)$. Then, the sum $\bar{\om}_t(G):=\sum_{k=0}^\infty \gamma_t^{k}\,\bar{\om}_t^{(k)}(\tilde{g}^k_t)$ is well defined, and when $t\to \infty$,
\begin{equation}
\bar{\om}_t(G) = \sum_{k=0}^\infty \gamma_t^k\,\bar{\om}_t^{(k)}(\tilde{g}^k_t) \cvlaw \sum_{k=0}^\infty \int_{\Delta_k}\int_{\mathbb{R}^k} g^k(\mathbf{t},\mathbf{x}) \eta^{\otimes k}(\mathrm{d}\mathbf{t},\mathrm{d}\mathbf{x}) = I(G).
\end{equation}
The convergence can be extended to a joint convergence for any finite collection $G_1, \dots , G_m$ in $ \bigoplus_{k\geq 0} L^2(\Delta_k\times \mathbb{R}^k)$:
\begin{equation} \label{eq:cor_ustat_jointCV}
\big(\bar{\om}_t(G_1), \dots ,\bar{\om}_t(G_m) \big) \cvlaw \big(I(G_1), \dots ,I(G_m) \big).
\end{equation}
\end{corollary}
\begin{remark}
Note that all functions in the corollary are defined on simplexes, as we have defined the Fock space for such functions.
\end{remark}

We first state a lemma that we will use several times. For a proof of the result, see \cite[Ch. 1. Th. 3.2]{billing}.

\begin{lemma} \label{lemme_billy}
Let $(S,\mathcal{S})$ be a metric space with his Borel $\sigma$-field. Suppose that $(X_t^n,Y_t)$ for $t\geq 0$, $n\in \mathbb{N}$ are random variables on $S^2$ and assume that the following diagram holds :
\[\begin{CD} 
X_t^n @>{(d)}>{t \to \infty}> Y^n\\
@V{\mathbb{P}, \text{ unif in } t}V{n\to \infty}V @V{(d)}V{n \to \infty}V\\
Y_t @. Y
\end{CD}\]
then $Y_t \cvlaw Y$.
\end{lemma}

\begin{proof}[Proof of the theorem]
We follow \cite{janson_1997}. In particular, we will focus in the first place on $k=1$ and $g\in L^1 \cap L^2([0,1]\times \mathbb{R})$.

\paragraph{$\mathbf{k=1}$ case.}
Let $g\in L^1 \cap L^2([0,1]\times \mathbb{R})$. When $k=1$, we have $\bar{\om}^{(1)}=\bar{\om}$, so we can use the complex form of the exponential formula for Poisson point processes (see equation \eqref{eq:PoissonExpMoment} below) to compute the characteristic function of $\gamma_t \,\bar{\om}_t(\tilde{g}_t)$. Let $u\in \mathbb{R}$, we have
\begin{align*}&\IP\left[e^{iu\gamma_t\bar{\om}_t(\tilde{g}_t)}\right]\\
&=\exp\left(\int_{[0,t]} \int_\mathbb{R}  \big(e^{iu\gamma_t g(s/t,x/\sqrt{t})}-1 -iu\gamma_t \, g(s/t,x/\sqrt{t})\, \big) \nu_t \, \mathrm{d}s \mathrm{d}x \right)\\
&=\exp\left(\int_{[0,1]} \int_\mathbb{R} \nu_t t^{3/2} \big(e^{iu\gamma_t g(s,x)}-1 -iu\gamma_t \, g(s,x)\big)\mathrm{d}s \mathrm{d}x \right),
\end{align*}
where the last equality comes from a change of variables.

By Taylor-Lagrange formula, we obtain
\begin{equation*}
\nu_t t^{3/2}\left|e^{iu\gamma_t g(s,x)}-1 -iu\gamma_t \, g(s,x)\right|\leq \nu_t t^{3/2} \gamma_t^2 \frac{u^2}{2} g(s,x)^2,
\end{equation*}
which gives $L^1$ domination, since $g$ is square-integrable and since  conditions (a) and (b) imply that $ \nu_t \gamma_t^2 \sim t^{-3/2}$.

Using again this asymptotic equivalence and the fact that $\gamma_t \to 0$, we get that the integrand converges pointwise to the function $(s,x) \mapsto -\frac{u^2}{2}g^2(s,x)$. Therefore, dominated convergence proves that, as $t\to\infty$,
\[\IP\left[e^{iu\gamma_t\tilde{\om}_t(\tilde{g}_t)}\right] \to \exp\left(-\frac{u^2}{2} \Vert g\Vert^2_2\right).\]

The limiting term is the Fourier transform of a centered Gaussian random variable of variance $\Vert g\Vert^2_2$, which has the same law as $I_1(g)$. This proves the first part of the theorem in the $k=1$ and $L^1\cap L^2$ case.

To prove the second part, we use the Cram\'er-Wold device which tells us that for a collection of real random variables, it is equivalent to show convergence in distribution of all finite linear combinations or to show joint convergence.

Thus, let $\alpha_1, \dots , \alpha_m \in \mathbb{R}$ and $g^1,\dots ,g^m \in L^1\cap L^2([0,1]\times \mathbb{R})$. By linearity of the stochastic integrals
\begin{align*}
\sum_{i=1}^m \alpha_i \, \gamma_t \, \bar{\om}_t (\tilde{g}^i_t) = \gamma_t \  \bar{\om}_t \left( \sum_{i=1}^m \alpha_i \, \tilde{g}^i_t \right) 
\cvlaw I_1\left(\sum_{i=1}^m \alpha_i \, g^i\right) = \sum_{i=1}^m \alpha_i \, I_1(g^i),
\end{align*}
where the convergence as $t\to \infty$ is ensured by the foregoing, since the combination $\sum_{i=1}^m \alpha_i g^i$ is square-integrable. By the Cram\'er-Wold device and as $t\to \infty$,
\[\big(\gamma_t\,\bar{\om}_t(\tilde{g}^{\,1}_t),\dots ,\gamma_t \,\bar{\om}_t(\tilde{g}^{\,m}_t)\big) \cvlaw \big(I_1(g^1),\dots ,I_1(g^m)\big).\]

\paragraph{$\mathbf{k\geq 1}$ case.}
Let $k\geq 1$ and let $g^1, \dots , g^k$ be the indicator functions of disjoint, finite and measurable sets $A_1, \dots , A_k \subset [0,1]\times \mathbb{R}$, and consider
\begin{equation} \label{eq:forme_prod}
g(\mathbf{s},\mathbf{x}) = g^1(s_1,x_1)\dots g^k(s_k,x_k).
\end{equation}

Equation \eqref{eq:formule_produit} writes \[\gamma_t^k\,\bar{\om}_t^{(k)}(\tilde{g}_t) = \prod_{i=1}^k \gamma_t\, \bar{\om}_t(\tilde{g}_t^i),\]
so joint convergence of the $\gamma_t \,\bar{\om}_t(\tilde{g}^i_t)$, $i\leq k$, from the $k=1$ case, implies that
\begin{equation*}
\bar{\om}_t^{(k)}(\tilde{g}_t) \cvlaw\prod_{i=1}^k I_1(g^i) = I_k(g),
\end{equation*}
where the equality comes from property \eqref{eq:IkOrthog} and the fact that the $g^i$'s are orthogonal in $L^2$. 

In fact, if one takes $g^1,\dots ,g^m$ of the form \eqref{eq:forme_prod}, so that they write $g^{i}(\mathbf{s},\mathbf x) =\prod_{j=1}^{k_i} g^{i,j}(s_j,x_j)$, the same argument, combined with the joint convergence of the $\gamma_t\,\bar{\om}_t(\tilde{g}^{\,i,j}_t)$ for  $1\leq i\leq m$ and $1\leq j\leq k_i$, proves that
\begin{equation} \label{cvloi_Vk}
\Big(\gamma_t^{k_1}\, \bar{\om}_t^{(k_1)}(\tilde{g}_t^{\,1}),  \dots  ,\gamma_t^{k_m}\,\bar{\om}_t^{(k_m)}(\tilde{g}_t^{\,m}) \Big)\\
\cvlaw \big(I_{k_1}(g^1),\dots ,I_{k_m}(g^m)\big).
\end{equation}

Now, denote by $V_k$ the linear subspace of $L^2([0,1]^k\times\mathbb{R})$ spanned by the functions of the form \eqref{eq:forme_prod}, with fixed dimension $k$. By linear combinations, \eqref{cvloi_Vk} can be extended for any collection $\left(g^i \in V_{k_i}\right)_{1\leq i\leq m}$, so this proves the whole theorem for functions of $V_k$, $k\geq 1$.

It is a standard result that $V_k$ is dense in $L^2([0,1]^k\times\mathbb{R}^k)$ for all $k\geq 1$. Let then $g$ be any function of $L^2([0,1]^k\times\mathbb{R}^k)$ and $(g^n)_{n\geq 1}$ be a sequence of functions of $V_k$ that converges to $g$ in $L^2$ norm. Conditions (a) and (b) imply that $\nu_t \, t^{3/2} \gamma_t^{2} \sim 1$. Hence, by the covariance structures and the linearity of $\bar{\om}_t^{(k)}$, we obtain for large enough $t$:
\begin{align*}
\IP\left[ \big(\gamma_t^k \bar{\om}_t^{(k)}(\tilde{g}^n_t)-\gamma_t^k\bar{\om}_t^{(k)}(\tilde{g}_t)\big)^2\right] &= \nu_t^k \gamma_t^{2k} k! \Vert (g^n-g)(\cdot/t,\cdot/\sqrt{t})\Vert^2_{L^2([0,t]^k\times \mathbb{R}^k)} \\
&= t^{3k/2}\nu_t^k \gamma_t^{2k} k! \Vert g^n-g\Vert_{L^2([0,1]^k\times \mathbb{R}^k)}^2 \\
&\leq 2 k! \Vert g^n-g\Vert_{L^2([0,1]^k\times \mathbb{R}^k)}^2 \to 0,
\end{align*}
as $n\to\infty$. Similarly for $I_k$:
\[
\mathbb Q\left[\big( I_k(g^n) - I_k(g) \big)^2\right] = k! \Vert g^n-g\Vert_{L^2([0,1]^k\times \mathbb{R}^k)}^2 \to 0
\]

We get the following diagram:
\[\begin{CD} 
\gamma_t \, \bar{\om}_t^{(k)}(\tilde{g}^n_t) @>{(d)}>{t \to \infty}> I_k(g^n)\\
@V{L^2, \text{ unif in } t}V{n\to \infty}V @V{(d)}V{n \to \infty}V\\
\gamma_t \, \bar{\om}_t^{(k)}(\tilde{g}_t) @. I_k(g)
\end{CD}\]
so that $\gamma_t \, \bar{\om}_t^{(k)}(\tilde{g}_t) \to I_k(g)$ by Lemma \ref{lemme_billy}. This proves the first part of the theorem. The joint convergence can be shown by the same argument, using again the Cram\'er-Wold device and approaching any linear combinations of Wiener-It\^o integrals of $L^2$ functions with linear combinations of integrals of $V_k$ functions.
\end{proof}

\begin{proof}[Proof of the corollary]
We focus on the first part of the corollary, since the joint convergence follows from the Cram\'er-Wold device and linearity of $\bar{\om}_t(G)$ and $I(G)$.

First and by definition, we know that $\sum_{k=0}^M I_k(g^k) \cvLdeux \sum_{k=0}^\infty I_k(g^k)$ as $M\to \infty$. Moreover, since we are now dealing with functions on the simplex, equation \eqref{eq:iso_simpl} leads to
\begin{equation*} \Vert \bar{\om}_t^{(k)}(\tilde{g}^k_t)\Vert_2^2 = \nu_t^k\Vert g^k(\cdot/t,\cdot/\sqrt{t})\Vert^2_{L^2(\Delta_k(0,t)\times \mathbb{R}^k)}  = \nu_t^k \, t^{3k/2} \Vert g^k\Vert^2_{L^2(\Delta_k\times \mathbb{R}^k)}.
\end{equation*} Conditions (a) and (b) imply that $\nu_t \, t^{3/2} \sim \gamma_t^{-2}$. Hence, as $\Vert g^k\Vert_2^2$ is summable, we obtain by absolute convergence that, uniformly in $t$ and as $M\to \infty$,
\[\sum_{k=0}^M \gamma_t^k \, \bar{\om}_t^{(k)}(\tilde{g}^k_t) \cvLdeux \sum_{k=0}^\infty \gamma_t^k \, \bar{\om}_t^{(k)}(\tilde{g}^k_t).\]

Moreover, it is a consequence the joint convergence part of the theorem that, for all $M$ and when $t\to\infty$,
\[\sum_{k=0}^M \gamma_t^k \, \bar{\om}_t^{(k)}(\tilde{g}^k_t) \cvlaw \sum_{k=0}^M I_k(g^k).\] 
Putting things together, we get the following diagram:
\[\begin{CD} 
\sum_{k=0}^M \gamma_t^k \, \bar{\om}_t^{(k)}(\tilde{g}^k_t) @>{(d)}>{t \to \infty}> \sum_{k=0}^M I_k(g^k)\\
@V{L^2, \text{ unif in } t}V{M\to \infty}V @V{(d)}V{M \to \infty}V\\
\sum_{k=0}^\infty \gamma_t^k \, \bar{\om}_t^{(k)}(\tilde{g}^k_t) @. \sum_{k=0}^\infty I_k(g^k),
\end{CD}\]
so by Lemma \ref{lemme_billy},
\[\bar{\om}_t(G) = \sum_{k=0}^\infty \gamma_t^k \, \bar{\om}_t^{(k)}(\tilde{g}^k_t) \cvlaw \sum_{k=0}^\infty I_k(g^k).
\]
\end{proof}

\section{Proofs}
\label{sec:proofs}
\subsection{Some useful formulas}
\begin{itemize}
\item 
For all non-negative and all non-positive measurable functions $h$, defined on $\mathbb{R}_+\times \mathbb{R}^d$, the Poisson formula for exponential moments (chapter 3. of \cite{LastPenrose}) writes
\begin{equation} \label{eq:PoissonExpMoment}
\P\left[ e^{\int h(s,x) \om_t(ds dx)}\right] = \exp \int_{]0,t]\times \mathbb{R}} \nu \rmd s \rmd x \left( e^{h(s,x)}-1 \right)   \;.
\end{equation}
The formula remains true when $h$ is replaced by $i h$, for any real integrable function $h$.
\item For all $s\geq 0$, we have
\begin{equation} \label{eq:formula_int_chi}
\int_\mathbb{R} \chi^r_{s,x} \hspace{0.3mm} \rmd x = r.
\end{equation}
\end{itemize}
\subsection{Proof of Theorem \ref{th:she} : SHE in the Poisson setting}
\label{sec:proofShe}
\begin{proof} Let $\xi_t = \exp(\b \om(V_t(B)) - \lambda(\b)\nu r^d t)$ and observe that 
\begin{equation*}\int_\mathbb{R} W(t,x)\varphi(x)\rmd x = \DP[\xi_t \varphi(B_t)].
\end{equation*} 
Then, recalling that $\om(V_t(B)) = \int \chi_{s,x} \, \om_t(\rmd s \rmd x)$, we use It\^o's formula \cite[Section II.5.]{IkedaWatanabe} to get that
\begin{align}
\nonumber \xi_t &= 1 -\lambda \nu r^d \int_0^t \xi_s \rmd s + \lambda \int_{(0,t]\times \mathbb{R}} \xi_{s-} \chi_{s,x} \, \om(\rmd s \rmd x)\\
& = 1 + \lambda\int_{(0,t]\times \mathbb{R}} \xi_{s-} \chi_{s,x} \, \bar{\om}(\rmd s \rmd x),
\end{align} as almost surely, $\IP$-a.s. $\xi_s = \xi_{s-}$ a.e.

As a difference of two increasing processes, $\xi$ is of finite variation over all bounded time intervals. Also note that one can get an expression to the measure associated to $\xi$ from the last equation. By the integration by part formula \cite[p.52]{jacod2013limit},
\[ \xi_t \varphi(B_t) = \xi_0 \varphi(B_0) + \int_0^t \xi_{s-} \rmd \varphi(B_s) + \int_0^t \varphi(B_s) \rmd \xi_s + [\xi,\varphi(B)]_t,\]
where $[\xi,\varphi(B)]_t=0$ since $\varphi(B)$ is continuous. Applying It\^o 's formula on $\rmd \varphi(B)$ and then taking $\DP$-expectation (which cancels the martingale term in the It\^o formula), one obtains that $\IP$-a.s.
\begin{align*}
&\int_\mathbb{R} W(t,x)\varphi(x)\rmd x\\
&  = \varphi(0) + \frac{1}{2}\int_0^t \DP[\xi_{s-} \Delta \varphi(B_s)] \rmd s + \lambda \int_{(0,t]\times\mathbb R} \DP[\varphi(B_s) \xi_{s-} \chi_{s,y}]\bar{\om}(\rmd s \rmd y)\\
& = \varphi(0) + \frac{1}{2}\int_0^t \int_\mathbb R \Delta \varphi(x) W(s-,x) \rmd x \rmd s \\
&  \qquad \qquad \qquad + \lambda \int_{(0,t]\times\mathbb R} \left(\int_\mathbb R \varphi(x) \mathbf{1}_{|y-x|\leq r/2} W(s-,x) \rmd x\right) \bar{\om}(\rmd s \rmd y).
\end{align*}
To conclude the proof, observe that we can apply Fubini's theorem to the last integral since for all $t>0$,
\begin{align*}
\IP \int_{(0,t] \times \mathbb{R}} \DP[|\varphi(B_s)|\xi_{s-} \chi_{s,y}] \om(\rmd s \rmd y) &= \nu e^\beta \int_{(0,t] \times \mathbb{R}} \IP[\xi_{s-}] \DP[|\varphi(B_s)| \chi_{s,y}] \rmd s \rmd y\\
& = \nu e^\beta r \int_0^t \DP[|\varphi(B_s)|] \rmd s < \infty,
\end{align*}
where we have used the Mecke equation \cite[4.1]{LastPenrose} in the first equality.
\end{proof}

\subsection{Proof of Theorem \ref{th:cv_Wt_interReg} : convergence of the P2L partition function}
\label{sec:PartitionFunctionsCV}
Introduce the following time-depending functions of $[0,1]^k\times\mathbb{R}^k$:
\begin{equation}
\phi_t^k(\mathbf{s},\mathbf{x}) = \gamma_t^{-k} \, \lambda(\b_t) ^k \DP\left[\prod_{i=1}^k \chi_{s_i,x_i}^{r_t/\sqrt{t}} (B)\right] \mathbf{1}_{\Delta_k}(\mathbf{s},\mathbf{x}).
\end{equation} 
Note that for all $(s,x)$, the diffusive scaling property of the Brownian motion implies that
\[\chi_{s/t,x/\sqrt{t}}^{r_t/\sqrt{t}} = \mathbf{1}_{|B_{s/t} - x/\sqrt{t}| \leq r_t/2\sqrt{t}} \eqlaw \chi^{r_t}_{s,x}.\]
Therefore, using notation $\widetilde{\phi}_t^k := (\widetilde{\phi_t^k})_t = \phi_t^k(\cdot/t,\cdot/\sqrt{t})$, we see that after simple rescaling, equation \eqref{eq:TkWtExpression} gives
\begin{equation} \label{eq:phitktilde} \gamma_t^k \, \widetilde{\phi}_t^k = T_k W_t \mathbf{1}_{\Delta_k(0,t)}.
\end{equation}
Besides, observe that by the symmetric property of $T_k W_t$ and the invariance of the Wiener-It\^o integrals under symmetrization \eqref{eq:sym_om}, we obtain\footnote{Note that from now on, we will always assume that $\om\eqlaw\om^{\nu_t}$, even if we drop the superscript notation.}:
\[\om_t^{(k)} \big(T_k W_t\big) = k! \, \om_t^{(k)} \big(T_k W_t \mathbf{1}_{\Delta_k(0,t)}\big).\]
Hence, Proposition \eqref{prop:Poisson_chaos_decomposition} and equation \eqref{eq:phitktilde} lead to the following expression of $W_t$:
\begin{equation} \label{eq:WtExpansionTwo}
W_t = \sum_{k=0}^\infty \gamma_t^k \, \bar{\om}^{(k)}_t \left(\widetilde{\phi}_t^k\right).
\end{equation}
Considering from now on $\phi_t^k$ as a function of the simplex, this writing of $W_t$ is of the type $\bar{\om}_t(G)$ (cf. Corollary \ref{coro_ustat}), although there is a time dependence in the $\phi_t^k$ functions. The purpose of the two following lemmas is to study the asymptotic behavior of these functions, as $t\to\infty$.
\subsubsection{Approximations in $L^2$-norm}
\label{sec:apprL2}
\begin{lemma} \label{lemma_L2} Let $k$ be a positive integer. We have the following properties:
\begin{enumerate}[label= (\roman*)]
\item For all $\mathbf{s} \in \Delta_k$, there exists a non-negative function $h_\mathbf{s} \in L^2(\mathbb{R}^k)$, such that 
\[ \forall \varepsilon \in  (0,1], \  \mathbf{x} \in \mathbb{R}^k, \quad \varepsilon^{-k}\,\mathrm{P}\left[\prod_{i=1}^k\chi_{s_i,x_i}^{\varepsilon}(B)\right] \leq  h_\mathbf{s}(\mathbf{x}).\]
\item There exists a non-negative function $H\in L^2(\Delta_k)$, such that
\[\forall \varepsilon >0, \  \forall \mathbf{s} \in \Delta_k, \quad \int_{\mathbb{R}^k} \left(\varepsilon^{-k}\, \mathrm{P}\bigg[\prod_{i=1}^k\chi_{s_i,x_i}^{\varepsilon}(B)\bigg] \right)^2 \mathrm{d}\mathbf{x} \leq H(\mathbf{s}).\]
\item \label{item:lemma_approx} We have the pointwise convergence, as $\varepsilon \to 0$, 
\[\forall \mathbf{s} \in \Delta_k, \, \forall \mathbf{x} \in \mathbb{R}^k, \quad \varepsilon^{-k} \, \mathrm{P}\left[\prod_{i=1}^k\chi_{s_i,x_i}^{\varepsilon}(B)\right] \to \rho^k(\mathbf{s},\mathbf{x}).\]
\end{enumerate}
\end{lemma}

\begin{proof} We use the convention $s_0=y_0=x_0=u_0$. 

(i) By Markov property of the Brownian Motion,
\begin{align} 
&\varepsilon^{-k} \, \mathrm{P}\left[\prod_{i=1}^k\chi_{s_i,x_i}^{\varepsilon}(B)\right] \nonumber\\
& =  \,\varepsilon^{-k}\int \limits_{\mathbb{R}^k} \prod_{i=1}^k \mathbf{1}_{|x_i-y_i|\leq \varepsilon/2} \, \rho(s_i - s_{i-1},y_i - y_{i-1}) \mathrm{d}\mathbf{y} \label{dvpt_markov} \\
&= \int\limits_{[-\frac{1}{2},\frac{1}{2}]^k} \prod_{i=1}^k \rho\big(s_i - s_{i-1},x_i - x_{i-1} + \varepsilon(u_i-u_{i-1}) \big) \mathrm{d}\mathbf{u},\label{dvpt_markov2}
\end{align}
where we have taken $y_i = x_i + \varepsilon u_i$.

Observe that, for $0<\varepsilon \leq 1$ and $u \in [-1,1]$,
\begin{equation*}
\frac{e^{-\frac{(x+\varepsilon u)^2}{2s}}}{\sqrt{2\pi s}}=\rho(s,x)e^{-\varepsilon\frac{2xu}{2s}}e^{-\varepsilon^2 \frac{u^2}{2s}}
\leq \rho(s,x) e^{|x|/s},
\end{equation*}

which leads to the following domination : \begin{equation} \label{ineg_lm}
\varepsilon^{-k} \, \mathrm{P}\left[\prod_{i=1}^k\chi_{s_i,x_i}^{\varepsilon}(B)\right]\leq \prod_{i=1}^k \rho(s_i - s_{i-1},x_i - x_{i-1})e^{|x_i - x_{i-1}|/(s_i - s_{i-1})}.
\end{equation} 

Define $h_\mathbf{s}(\mathbf{x})$ to be the right-hand side of \eqref{ineg_lm}, so that what is left to prove is that $h_\mathbf{s} \in L^2(\mathbb{R}^k)$ for $\mathbf{s}\in \Delta_k$. With the change of variables $z_i = x_i - x_{i-1}$ of Jacobian $J=1$ and by Tonelli's theorem
\begin{align*} 
\int_{\mathbb{R}^k} h_\mathbf{s}(\mathbf{x})^2 \mathrm{d}\mathbf{x} &= \int_{\mathbb{R}^k}\prod_{i=1}^k \frac{e^{-(x_i-x_{i-1})^2/(s_i - s_{i-1})}}{2\pi (s_i - s_{i-1})}e^{2|x_i-x_{i-1}|/{(s_i - s_{i-1})}} \mathrm{d}\mathbf{x} \\
&= \prod_{i=1}^k\int_{\mathbb{R}} \frac{e^{-z_i^2/(s_i - s_{i-1})}}{2\pi (s_i - s_{i-1})}e^{2|z_i|/{(s_i - s_{i-1})}} \mathrm{d}z_i,
\end{align*}
which is finite as each integral converges.

(ii) One can first note that for all $s>0$, $\rho(s,x) \leq 1/{\sqrt{2\pi s}}.$
This combined with equation \eqref{dvpt_markov}  gives us, for all $\mathbf{s}\in \Delta_k$,
\begin{align}\mathrm{P}\left[\prod_{i=1}^k\chi_{s_i,x_i}^{\varepsilon}(B)\right] &\leq (2\pi)^{-k/2} \prod_{i=1}^k \frac{1}{\sqrt{s_i - s_{i-1}}} \int_{\mathbb{R}^k} \prod_{i=1}^k \mathbf{1}_{|x_i-y_i|\leq \varepsilon/2} \, \mathrm{d}\mathbf{y} \nonumber\\
& \leq \varepsilon^{k} \prod_{i=1}^k \frac{1}{\sqrt{s_i - s_{i-1}}}.
\end{align}
Let $H(\mathbf s)=\prod_{i=1}^k (s_i - s_{i-1})^{-1/2}$ be the product appearing in the last inequality. We saw in the proof of Proposition \ref{prop_rhoIsInFock} that $H$ is an element of $L^1(\Delta_k)$. Furthermore:
\begin{align}
\int_{\mathbb{R}^k} \left(\varepsilon^{-k}\, \mathrm{P}\bigg[\prod_{i=1}^k\chi_{s_i,x_i}^{\varepsilon}(B)\bigg] \right)^2 \mathrm{d}\mathbf{x}
\leq & H(\mathbf s) \int_{\mathbb{R}^k} \varepsilon^{-k} \, \mathrm{P}\bigg[\prod_{i=1}^k\chi_{s_i,x_i}^{\varepsilon}(B)\bigg] \mathrm d \mathbf x \nonumber \\
\label{ineq_lemma} = & H(\mathbf s),
\end{align}
where we have used Tonelli's theorem in the equality.

(iii) This result can be derived from equation \eqref{dvpt_markov2}, using continuity in $x$ of $\rho(s,x)$ for a fixed $s>0$.
\end{proof}

From the last lemma, we can derive $L^2$ properties of $\phi^k_t$:
\begin{lemma} \label{lemma_cvL2} Let $k$ be a positive integer. We have:
\begin{enumerate}[label= (\roman*)]
\item The following convergence holds: \begin{equation*}\Vert\phi^k_t - (\b^*)^k \rho^k\Vert_{L^2(\Delta_k\times \mathbb{R}^k)} \underset{t \to \infty}{\longrightarrow} 0.
\end{equation*}
\item There exists a positive constant $C=C(\b^*)$, such that
\begin{equation*}
\sup_{t\in [0,1]}\Vert\phi^k_t\Vert_{L^2([0,1]^k\times\mathbb{R}^k)} \leq C^k\Vert\rho^k\Vert_{L^2(\Delta_k\times\mathbb{R}^k)}.
\end{equation*}
\end{enumerate}
\end{lemma} 
\begin{proof} (i)
Recall that on $\Delta_k(0,t) \times \mathbb{R}^k$:
\begin{equation*}\phi_t^k(\mathbf{s},\mathbf{x})=\gamma_t^{-k} \, \lambda ^k \DP\left[\prod_{i=1}^k \chi_{s_i,x_i}^{r_t/\sqrt{t}} (B)\right].
\end{equation*}
Conditions (a) and (b) imply that, as $t\to\infty$,
\begin{equation}\label{eq:lemma_equivalent}\gamma_t^{-1} \lambda(\b_t) \sim \b^* \frac{\sqrt{t}}{r_t},
\end{equation} which leads to the existence of a constant $c>1$, such that, for $t$ large enough and all $k\geq 1$:
\begin{equation} \label{ineq_lemma2} \gamma_t^{-k} \lambda^k \leq c^k |\b^*|^k \left(\frac{r_t}{\sqrt{t}}\right)^{-k}.
\end{equation}
Observe that $r_t/\sqrt{t} \to 0$ by condition (c). Then, Lemma \ref{lemma_L2} implies the existence of two non-negative dominating functions $h_\mathbf{s} \in L^2(\mathbb{R}^k)$ and $H \in L^1(\Delta_k)$, such that, for large enough $t$ and $\mathbf{s}\in \Delta_k$,
\begin{equation} \label{first_domi}
\forall \mathbf{x} \in \mathbb{R}^k \quad {\phi_t^k}(\mathbf{s},\mathbf{x}) \leq h_\mathbf{s}(\mathbf{x}),
\end{equation}
and
\begin{equation} \label{sec_domi}
\int_{\mathbb{R}^k} {\phi_t^k}(\mathbf{s},\mathbf{x})^2 \mathrm{d}\mathbf{x} \leq H(\mathbf{s}).
\end{equation}
Furthermore, because of the equivalence \eqref{eq:lemma_equivalent}, point \ref{item:lemma_approx} of the same lemma shows that we have the pointwise convergence:
\[\forall \mathbf{s} \in [0,1], \, \forall \mathbf{x} \in \mathbb{R}^k, \quad \phi_k^t(\mathbf{s},\mathbf{x}) \underset{t\to \infty}{\longrightarrow} (\b^*)^k \rho^k(\mathbf{s},\mathbf{x}).\]

From \eqref{first_domi}, we get by the dominated convergence theorem that
\begin{equation*} \forall \mathbf s \in \Delta_k, \quad \int_{\mathbb{R}^k} \big({\phi_t^k}(\mathbf{s},\mathbf{x})-(\b^*)^k\rho^k(\mathbf{s},\mathbf{x})\big)^2 \mathrm{d}\mathbf{x} \underset{t\to \infty}{\longrightarrow} 0,
\end{equation*}and as we have
\[\int_{\mathbb{R}^k}\big({{\phi_t^k}}(\mathbf{s},\mathbf{x})-(\b^*)^k\rho^k(\mathbf{s},\mathbf{x})\big)^2 \mathrm{d}\mathbf{x}\leq \int_{\mathbb{R}^k} 2{{\phi_t^k}}(\mathbf{s},\mathbf{x})^2+2(\b^*)^{2k}\rho^k(\mathbf{s},\mathbf{x})^2 \mathrm{d}\mathbf{x},\] 
where the right hand side is dominated in $L^1(\Delta_k)$ using equation \eqref{sec_domi}, we can use the dominated convergence theorem and obtain the convergence
\begin{equation*}
\Vert\phi_k^t - (\b^*)^k\rho^k\Vert^2_{L^2(\Delta_k\times \mathbb{R}^k)}=\int_{\Delta_k} \left(\int_{\mathbb{R}^k}\big({\phi_k^t}(\mathbf{s},\mathbf{x})-(\b^*)^k\rho^k(\mathbf{s},\mathbf{x})\big)^2 \mathrm{d}\mathbf{x}\right) \mathrm d \mathbf{s} \underset{t\to \infty}{\longrightarrow} 0.
\end{equation*}

(ii) Using inequalities \eqref{ineq_lemma2} and  \eqref{ineq_lemma}, we get
\[ \Vert\phi_k^t\Vert_{L^2([0,1]^k\times \mathbb{R}^k)}^2 \leq c^k |\b^*|^k \int_{\Delta_k} \prod_{i=1}^k \frac{1}{\sqrt{s_i - s_{i-1}}}\mathrm d \mathbf s=C^{k} \Vert\rho^k\Vert_{L^2(\Delta_k\times \mathbb{R}^k)}^2,\]
with $C=2\sqrt{\pi}c|\b^*|$.
\end{proof}
\subsubsection{Proof of the theorem}
We are now ready to prove Theorem \ref{th:cv_Wt_interReg}. From Proposition \ref{prop_rhoIsInFock}, we know that $R(\b^*)=((\beta^*)^k\rho^k)_{k\geq 0}\in\bigoplus_{k\geq 0} L^2(\Delta_k\times \mathbb{R}^k)$, so we get from Corollary \ref{coro_ustat} that when $t\to \infty$,
\begin{equation} \sum_{k=0}^\infty (\beta^*)^k \gamma_t^k\,\bar{\om}_t^{(k)}(\tilde{\rho}^{\,k}_t) \cvlaw \mathcal{Z}_{\beta^*}.
\end{equation}

In addition to this, we saw at equation \eqref{eq:WtExpansionTwo} that the renormalized partition function writes
$
W_t = \sum_{k=0}^\infty \gamma_t^k \, \bar{\om}^{(k)}_t \left(\tilde{\phi}_t^k\right),
$
where $\phi^k_t \cvLdeux (\b^*)^k\rho^k$ by Lemma \ref{lemma_cvL2}. It is a standard result that if $X_n$ and $Y_n$ are real random variables such that $Y_n \cvlaw Y$ and $\Vert Y_n-X_n\Vert_2 \longrightarrow 0$, then $X_n \cvlaw Y$.
Hence, in order to prove that $W_t \cvlaw \mathcal{Z}_{\beta^*}$, it suffices to show that
\begin{equation}\left\Vert \sum_{k=0}^\infty \gamma_t^k \, \bar{\om}_t^{(k)}(\tilde{\phi}_t^k) - \sum_{k=0}^\infty (\b^*)^k \gamma_t^k\,\bar{\om}_t^{(k)}(\tilde{\rho}^{\,k}_t) \right\Vert_2^2 \underset{{t} \to {\infty}}{\longrightarrow} 0.
\end{equation}

By linearity of $\tilde{\om}_t^{(k)}$ and orthogonality of the terms for two different $k$, we get from Pythagoras' identity that the norm can be written as \[\sum_{k=0}^\infty \gamma_t^{2k} \, \Vert\bar{\om}_t^{(k)}\big(\phi_t^k(\cdot/t,\cdot/\sqrt{t})- (\b^*)^k\rho^{k}(\cdot/t,\cdot/\sqrt{t})\big)\Vert_2^2.\]

For all $g\in L^2(\Delta_k\times \mathbb{R}^k)$, equation \eqref{eq:iso_simpl} and a substitution of variables lead to
\begin{equation*}\Vert\bar{\om}_t^{(k)}\big(g(\cdot/t,\cdot/\sqrt{t})\Vert_2^2 = \nu_t^k\Vert g(\cdot/t,\cdot/\sqrt{t})\Vert^2_{L^2(\Delta_k(0,t)\times \mathbb{R}^k)}  =  \nu_t^k t^{3k/2} \Vert g\Vert^2_{L^2(\Delta_k\times \mathbb{R}^k)},
\end{equation*}
so that the above sum is given by
\[\sum_{k=0}^\infty \gamma_t^{2k} \nu_t^k t^{3k/2} \, \Vert \phi_t^k-(\b^*)^k\rho^{k}\Vert_{L^2(\Delta_k\times \mathbb{R}^k)}^2.\]

Conditions (a) and (b) imply that $\gamma_t^2 \nu_t^k t^{3/2} \sim 1$, so by lemma \ref{lemma_cvL2}, the summand tends to zero, as $t\to\infty$, and it is dominated by $C^{2k} \Vert\rho^k\Vert^2_2$, where $C=C(\b^*)$ is some positive constant. As this dominating sequence is summable (Proposition \ref{prop_rhoIsInFock}), the dominated convergence theorem concludes the proof.

\subsection{Proof of Theorem \ref{th:cv_P2PWt_interReg} : convergence of the point-to-point partition function}
Using again Theorem \ref{th:wienerItoChaos}, and after similar renormalization to what was done in the beginning of Section \ref{sec:PartitionFunctionsCV}, we find that
\begin{equation} \label{eq:P2PChaos} \sqrt{t} W(\b_t,tS,\sqrt{t} \hspace{0.3mm}Y;t\hspace{0.3mm}T,\sqrt{t}\hspace{0.3mm}X) = \sum_{k=0}^\infty \gamma_t^k \, \bar{\om}^{(k)}_t\left(\widetilde{\psi}_t^k(S,Y;T,X)\right),
\end{equation}
where
\begin{align*}
&\psi_t^k(S,Y;T,X)(\mathbf{s},\mathbf{x})\\ 
&=\gamma_t^{-k} \lambda(\b_t)^k \rho(T-S,X-Y) \DP_{S,Y}^{T,X} \left[\prod_{i=1}^k \chi_{s_i,x_i}^{r_t/\sqrt{t}} (B)\right] \mathbf{1}_{\Delta_k(S,T)}.
\end{align*}
Analogous calculations to those of Section \ref{sec:apprL2} will show that for all $k\geq 0$, as $t\to\infty$,
\begin{equation*}
\psi_t^k(S,Y;T,X) \cvLdeux {\b^*}^k \rho^k(\cdot\,;\,S,Y;T,X),
\end{equation*} 
where, by Corollary \ref{coro_ustat}, the limiting functions have the property that
\begin{equation} \label{eq:CVlawP2P}
\sum_{k=0}^\infty \gamma_t^k {\b^*}^k \, \bar{\om}^{(k)} \left(\tilde{\rho}^k_t(\cdot\,;\,S,Y;T,X)\right) \cvlaw \mathcal{Z}_{\b^*}(S,Y;T,X).
\end{equation}

The theorem then follows from similar arguments to those of the proof of the convergence of the point-to-line partition function, that is by
showing that the right-hand side of \eqref{eq:P2PChaos} and the left-hand side of \eqref{eq:CVlawP2P} are close in $L^2$ norm.

\subsection{Proof of Theorem \ref{th:cvlaw} : convergence in terms of processes}
\label{sec:proofCVProcesses}
In order to show tightness of $\mathcal{Y}_t$, we will rely on Mitoma's criterion \cite{mitoma1983tightness,Walsh}. It will help us reduce the problem of showing tightness of a two variables process, to the problem of showing tightness of a set of one variable processes.

In what follows, for any function $F\in D\left([0,1],\mathcal{D}'(\mathbb{R})\right)$ and $\varphi\in \mathcal{D}(\mathbb{R})$, we set
\begin{equation}
F(T,\varphi) := \int F(T,X)\varphi(X)\rmd X.
\end{equation}
\begin{proposition}[\textup{\cite{mitoma1983tightness}}] \label{th:Mitoma}
Let $(F_t)_{t\geq 0}$ be a family of processes in $D\left([0,1],\mathcal{D}'(\mathbb{R})\right)$. If, for all $\varphi \in \mathcal{D}(\mathbb{R})$, the family $T\to F_t(T,\varphi),t\geq 0$ is tight in the real cadl\`ag functions space $D([0,1],\mathbb{R})$, then $(F_t)_{t\geq 0}$ is tight in $D\left([0,1],\mathcal{D}'(\mathbb{R})\right)$.
\end{proposition}

Then, in order to prove uniqueness of the limit, we use the following proposition:
\begin{proposition}[\textup{\cite{mitoma1983tightness}}] \label{prop:mitomaIdentification}
Let $(F_t)_{t\geq 0}$ be a tight family of processes in the space $D\left([0,1],\mathcal{D}'(\mathbb{R})\right)$.  If there exists a process $F\in D\left([0,1],\mathcal{D}'(\mathbb{R})\right)$ such that, for all $n\geq 1$, $T_1,\dots,T_n\in [0,1]$ and $\varphi_1,\dots,\varphi_n \in \mathcal{D}(\mathbb R)$, we have as $t\to\infty$:
\[\left(F_t(T_1,\varphi_1),\dots,F_t(T_n,\varphi_n)\right) \cvlaw \left(F_t(T_1,\varphi_1),\dots,F_t(T_n,\varphi_n)\right),
\]
then $F_t \cvlaw F$.
\end{proposition}

\subsubsection{Identification of the limit}
\begin{proposition} Let $\varphi\in\mathcal{D}(\mathbb R)$. Then, for all $T\geq 0$ and as $t\to\infty$,
\begin{equation} \label{eq:identifOfTheLimit}
\mathcal{Y}_t(T,\varphi) := \int \mathcal{Y}_t(T,X) \varphi(X) \rmd X \cvlaw \int \mathcal{Z}_{\beta^*}(T,X) \varphi(X) \rmd X.
\end{equation}
Moreover, the convergence extends to a joint convergence as in Proposition \ref{prop:mitomaIdentification}.
\begin{proof}
Once again, we rely on Theorem \ref{th:wienerItoChaos} and similar renormalization to the beginning of Section \ref{sec:PartitionFunctionsCV} to get that
\begin{equation*} \mathcal{Y}_t(T,\varphi) = \sum_{k=0}^\infty \gamma_t^k \, \bar{\om}^{(k)}\left(\widetilde{\psi}_t^k(T,\varphi)\right),
\end{equation*}
where
\begin{align*}
& \psi_t^k(T,\varphi)(\mathbf{s},\mathbf{x})\\
&= \gamma_t^{-k} \lambda(\b_t)^k \int_\mathbb R \rho(T,X) \DP_{0,0}^{T,X} \left[\prod_{i=1}^k \chi_{s_i,x_i}^{r_t/\sqrt{t}} (B)\right] \varphi(X) \rmd X \, \mathbf{1}_{\Delta_k}(\mathbf{s},\mathbf{x}).
\end{align*}
Then, for all $k\geq 0$ and as $t\to\infty$, we have:
\begin{equation*}
\psi_t^k(T,\varphi) \cvLdeux g^k := {\b^*}^k \int_\mathbb R \rho^k(\cdot\,;0,0;\,T,X) \varphi(X)\rmd X.
\end{equation*}
To see this, apply Cauchy-Schwarz's inequality to obtain that
\begin{align*}
&\Vert\psi_t^k(T,\varphi) - g^k\Vert^2_{2}  \\
&= \int_{\Delta_k\times\mathbb{R}^{k}} \Bigg( \int_{\mathbb{R}} \varphi(X) \Bigg(\gamma_t^{-k} \lambda(\b_t)^k\rho(T,X) \DP_{0,0}^{T,X} \left[\prod_{i=1}^k \chi_{s_i,x_i}^{r_t/\sqrt{t}}\right] \\
& \hspace{5cm} - {\b^*}^k \rho^k(\mathbf{s},\mathbf{x},0,0;T,X) \Bigg)\rmd X \Bigg)^2 \rmd \mathbf{s}\rmd \mathbf{x}\\
& \leq \Vert \varphi \Vert^2_{2}  \int_{{\Delta_k\times\mathbb{R}^{k+1}}} \Bigg(\gamma_t^{-k} \lambda(\b_t)^k\rho(T,X) \DP_{0,0}^{T,X} \left[\prod_{i=1}^k \chi_{s_i,x_i}^{r_t/\sqrt{t}}\right]\\
& \hspace{5.5cm} - {\b^*}^k \rho^k(\mathbf{s},\mathbf{x};0,0;T,X) \Bigg)^2 \rmd \mathbf{s}\rmd \mathbf{x}\rmd X.
\end{align*}
By similar estimates to those obtained in Section \ref{sec:apprL2}, we get that this last integral goes to $0$. 

Equation \eqref{eq:sheContreTestfunction} and Corollary \ref{coro_ustat}  imply that
\begin{equation*}
\sum_{k=0}^\infty \gamma_t^k \, \bar{\om}^{(k)} \left(\tilde{g}^k_t\right) \cvlaw \int_\mathbb{R}\mathcal{Z}_{\b^*}(S,X)\varphi(X)\rmd X,
\end{equation*}
so that convergence \eqref{eq:identifOfTheLimit} follows from the same arguments we used for the convergence of the P2L functions. Finally, the joint convergence can be obtained using \eqref{eq:cor_ustat_jointCV} in Corollary \ref{coro_ustat}.
 
\end{proof}
\end{proposition}

\subsubsection{Tightness}
The process $\Lambda^t$ is a Poisson point process of intensity measure $t^{3/2} \nu_t\, \rmd S \rmd X$.
By simple rescaling of the Poisson stochastic heat equation \eqref{eq:sheForPoisson}, one can write that
\begin{equation}
\mathcal{Y}_t(T,\varphi) = \varphi(0) + A^t_T(\varphi) + M^t_T(\varphi),
\end{equation}
where
\begin{align*}
A^t_T(\varphi) & = \frac{1}{2}\int_0^T\int_\mathbb{R} \Delta \varphi(X)  \mathcal{Y}_t(S,X) \rmd S \rmd X, \\
M^t_T(\varphi) & = \lambda(\b_t) \int_\mathbb{R} \varphi(X) \left(\int_{(0,T]\times \mathbb R} \mathcal{Y}_t(S-,X) \mathbf{1}_{|X-Y|\leq r_t/2\sqrt{t}}\,\overline{\Lambda^t}(\rmd S \rmd Y)\right)\rmd X.
\end{align*}
We will show that both $A^t(\varphi)$ and $M^t(\varphi)$ are tight in $\mathcal{D}([0,1],\mathbb{R})$. If this is proven, then $\big(A^t(\varphi),M^t(\varphi)\big)$ is tight, hence $\mathcal{Y}_t(\cdot,\varphi)$ is tight.
\paragraph{Tightness of $A^t(\varphi)$:}
To prove tightness, we will use Kolmogorov's criterion \cite[Theorem 12.3]{billingsley1968convergence}. For this we need estimates on the moments on the variations of $A^t(\varphi)$, which will be derived through the next lemma:
\begin{lemma} \label{lem:boundOnMomentsOfwt}
Suppose that $(\b_t)_{t\geq 0}$ is bounded by above. Then, there exists a constant $C = C(\b^*,p)$ verifying $0<C<\infty$ and that for $t$ large enough and all $p>1$, $T>0$ and $X\in\mathbb R$,
\begin{equation} 
\IP[\mathcal{Y}_t(T,X)^p] \leq C \rho(T,X)^p.
\end{equation}
\end{lemma}
Suppose for a moment that the lemma is proven, and let $p\geq 2$ be an integer and $U\leq T$ in $[0,1]$. We have:
\begin{align*}
\IP\left[\left|A^t_T(\varphi)-A^t_U(\varphi)\right|^p\right] & \leq 2^{-p} \IP\left[\left(\int_{[U,T]\times\mathbb R} \left|\Delta \varphi(X)\right|  \mathcal{Y}_t(S,X) \rmd S \rmd X\right)^p\right]\\
& \leq 2^{-p} \Vert \Delta \varphi \Vert^{p}_\infty \int_{[U,T]^p\times\mathbb{R}^p}\IP\left[\prod_{i=1}^p \mathcal{Y}_t(S_i,X_i)\right] \rmd \mathbf{S}\rmd \mathbf{X}.
\end{align*}
By Lemma \ref{lem:boundOnMomentsOfwt}, the functions $\mathcal{Y}_t(S_i,X_i)/\rho(S_i,X_i)$ are bounded in $L^{p}$, so we can use the generalized H\"older inequality to bound the expectation of the product in the right-hand side. We get that there is constant $C=C(p)>0$ such that
\begin{align*}
\IP\left[ \left| \left\langle A^t(\varphi) \right\rangle_T - \left\langle A^t(\varphi) \right\rangle_U\right|^p \right] & \leq C  \Vert \Delta \varphi \Vert^{p}_\infty \int_{[U,T]^p\times\mathbb{R}^p}\prod_{i=1}^p \rho(S_i,X_i) \rmd \mathbf{S}\rmd \mathbf{X}\\
& = C  \Vert \Delta \varphi \Vert^{p}_\infty \left( \int_{[U,T]\times\mathbb{R}} \rho(S,X)\rmd S \rmd X \right)^p\\
&= C \Vert \Delta \varphi \Vert^{p}_\infty \; | T - U |^p.
\end{align*}
This shows that $A^t(\varphi)$ verifies the assumptions of Kolmogorov's criterion.
\begin{proof}[Proof of Lemma \ref{lem:boundOnMomentsOfwt}]
Let $p>1$, $T>0$ and $X\in\mathbb R$.
We have:
\[\mathcal{Y}_t(T,X)^p = \rho(T,X)^p \DP^{\otimes p} \left[\exp \int_{[0,tT]\times\mathbb{R}} \sum_{i=1}^p \b_t\chi^{r_t}_{s,x}(X^i)\, \om(\rmd s \rmd x) \right]e^{-ptT\lambda \nu_t r_t},\]
where $X^1,\dots,X^p$ are independent Brownian bridges from $(0,0)$ to $(tT,\sqrt{t}X)$. By the exponential formula \eqref{eq:PoissonExpMoment}, we get
\[\frac{\IP[\mathcal{Y}_t(T,X)^p]}{\rho(T,X)^p} = \DP^{\otimes p} \left[\exp \int_{[0,tT]\times\mathbb{R}} \left(e^{\sum_{i=1}^p \b_t\chi^{r_t}_{s,x}(X^i)}-1\right)\, \nu_t\rmd s \rmd x \right]e^{-ptT\lambda \nu_t r_t}.\]
Then, observe that
\begin{align*}
e^{\sum_{i=1}^p \b_t\chi^{r_t}_{s,x}(X^i)} &= \prod_{i=1}^p\left(1+\lambda \chi^{r_t}_{s,x}(X^i)\right)\\
& = 1 + \sum_{i=1}^p \lambda \chi^{r_t}_{s,x}(X^i) + \sum_{k=2}^p \lambda^k \sum_{p_1 < \dots < p_k \leq p} \prod_{i=1}^k \chi^{r_t}_{s,x}(X^{p_i}).
\end{align*}
Using equation \eqref{eq:formula_int_chi}, we are left with:
\begin{equation} \label{eq:normLpOfwtx}
\begin{split}
&\frac{\IP[\mathcal{Y}_t(T,X)^p]}{\rho(T,X)^p}=\\
&\hspace{0.9em}  \DP^{\otimes p} \left[\prod_{k=2}^p \prod_{p_1 < \dots < p_k \leq p}\exp \left( \nu_t\lambda(\b_t)^k \int_{[0,tT]\times\mathbb{R}}  \prod_{i=1}^k \chi^{r_t}_{s,x}(X^{p_i})\, \rmd s \rmd x \right) \right].
\end{split}
\end{equation}
We claim that for each $q>0$ and $k\geq 2$, there exists a constant $C=C(q,\b^*)>0$, such that for all $t$ large enough and all $T,X$,
\begin{equation} \label{eq:boundLpLemmawtx}
\DP^{\otimes p} \left[\exp \left( q\nu_t|\lambda(\b_t)|^k \int_{[0,tT]\times\mathbb{R}}  \prod_{i=1}^k \chi^{r_t}_{s,x}(X^{p_i})\, \rmd s \rmd x \right) \right] \leq C.
\end{equation}
If this is proven, then the generalized H\"older inequality implies that the right-hand side of \eqref{eq:normLpOfwtx} is bounded, which is the claim of the lemma.

First, notice that for all $k\geq 2$,
\begin{align*}
\int_0^{tT} \int_\mathbb{R} \prod_{i=1}^k \chi^{r_t}_{s,x}(X^{p_i})\rmd x \rmd s & \leq \int_0^{tT} \int_\mathbb{R} \chi^{r_t}_{s,x}(X^{p_1})\chi^{r_t}_{s,x}(X^{p_2})\rmd x \rmd s\\
& \leq \int_0^{tT} \int_\mathbb{R} \mathbf{1}_{|X_s^{p_1}-X_s^{p_2}|\leq r_t/2}\,\chi^{r_t}_{s,x}(X^{p_1})\rmd x \rmd s\\
& = t r_t \int_0^{T} \mathbf{1}_{|X_{tS}^{p_1}-X_{tS}^{p_2}|\leq r_t/2} \, \rmd S
\end{align*}
As $X_{tS}^{p_1}-X_{tS}^{p_2} \eqlaw \sqrt{2t} \tilde{X}_S$ , where $\tilde{X}$ is a Brownian bridge $(0,0)\to(T,0)$, we can bound the left hand side of \eqref{eq:boundLpLemmawtx} by
\begin{equation*}\sum_{m=0}^\infty \frac{1}{m!} (qt r_t\nu_t |\lambda(\b_t)|^k)^m \DP\left[\left(\int_0^T \mathbf{1}_{|\tilde{X}_S|\leq r_t/\sqrt{8t}}\,\rmd S\right)^m\right].
\end{equation*}
By symmetry, we have:
\begin{align*}
&\frac{1}{m!}\DP\left[\left(\int_0^T \mathbf{1}_{|\tilde{X}_S|\leq r_t/\sqrt{8t}}\,\rmd S\right)^m\right]\\
& = \int_{\Delta_m(0,T)} \DP\left[\prod_{i=1}^m\mathbf{1}_{|\tilde{X}_{S_i}|\leq r_t/\sqrt{8t}} \right] \rmd \mathbf{S}\\
& = \int_{\Delta_m(0,T)} \int_{[-\frac{r_t}{2\sqrt{2t}},\frac{r_t}{2\sqrt{2t}}]} \frac{\prod_{i=1}^{m+1} \rho(S_i-S_{i-1},X_{i-1}-X_i)}{\rho(T,0)}\rmd \mathbf{S} \rmd \mathbf{X},
\end{align*}
where $S_0=X_0=X_{m+1}=0$ and $S_{m+1}=T$. Using that $\rho(S,X)\leq \sqrt{2\pi S}^{-1}$, we get that
\begin{align*}
&\frac{1}{m!}\DP\left[\left(\int_0^T \mathbf{1}_{|\tilde{X}_S|\leq r_t/\sqrt{8t}}\,\rmd S\right)^m\right]\\
 & \leq (\sqrt{2\pi})^{-m} \sqrt{T} \left(\frac{r_t}{\sqrt{2t}}\right)^m \int_{\Delta_m(0,T)} \prod_{i=1}^{m+1} {\sqrt{S_i-S_{i-1}}}^{-1}\rmd \mathbf{S}\\
& \leq C^m r_t^{m} t^{-m/2} T^{m/2} \int_{\Delta_m(0,1)} {\sqrt{U_i-U_{i-1}}}^{-1}\rmd \mathbf{U}\\
& = C^m r_t^{m} t^{-m/2} T^{m/2} \frac{\sqrt{\pi}^m}{\Gamma(m/2)},
\end{align*}
where $C>0$ is some constant and where the value of the integral in the third equation was identified via the Dirichlet distribution.

Since $\b_t$ is assumed to be bounded by above, $\lambda(\b_t)$ is bounded. Then, as $k\geq 2$, the scaling relation (a) implies that there exists a constant $C_1=C_1(\b^*)>0$, such that $t^{1/2} r_t^2 \nu_t\lambda(\b_t)^k \leq C_1$. Added to the fact that $T\leq 1$, we get that there exists a finite constant $C_2>0$, depending only on $\b^*$ and $q$, such that for $t$ large enough:
\[\sum_{m=0}^\infty \frac{1}{m!} (qt r_t\nu_t |\lambda|^k)^m \DP\left[\left(\int_0^T \mathbf{1}_{|\tilde{X}_S|\leq r_t/\sqrt{8t}}\,\rmd S\right)^m\right] \leq \sum_{m=0}^\infty \frac{C_2^m}{\Gamma(m/2)} < \infty.\]
This proves \eqref{eq:boundLpLemmawtx}, which ends the proof of the lemma.

\end{proof}

\paragraph{Tightness of $M^t(\varphi)$:}
The process $T\mapsto M^t_T(\varphi)$ is a martingale (see equation \eqref{eq:mtvarphi} below) with respect to the filtration induced by $(\om_{tT})_{T\in[0,1]}$. It is therefore possible to rely on Aldous' criterion to show tightness:

\begin{theorem}[Aldous' criterion for martingales \textup{\cite[Chap. VI Theorem 4.13]{jacod2013limit}}]
Let $(N^t)_{t\geq 0}$ be a family of martingales in $ D([0,1],\mathbb R)$. Assume that:
\begin{enumerate}[label=(\roman*)]
\item The family $(N^t_0)_{t\geq 0}$ is tight.
\item The family of previsible brackets $(\langle N^t \rangle)_{t\geq 0}$ is tight in $\mathcal C([0,1],\mathbb R)$.
\end{enumerate}
Then, the $(N^t)_{t \geq 0}$ are tight in $ D([0,1],\mathbb R)$.
\end{theorem}
\noindent In our case, point (i) is immediately verified, as $M^t_0(\varphi) = 0$. To show that point (ii) holds, we use Kolmogorov's criteria. We have:
\begin{gather}
M^t_T(\varphi) = \int_{(0,T]\times \mathbb R} f(S,Y,\om) \overline{\Lambda^t}(\rmd S \rmd Y),\label{eq:mtvarphi}\\
\text{where} \quad f(S,Y,\om) = \lambda(\b_t) \int_{\mathbb R} \mathcal{Y}_t(S-,X) \varphi(X)  \mathbf{1}_{|X-Y|\leq r_t/2\sqrt t} \, \rmd X.\nonumber
\end{gather}
Since $f$ is predictable and $\Lambda^t$ has intensity $t^{3/2} \nu_t \, \rmd S \rmd X$, the bracket can be expressed \cite[Section II.3.]{IkedaWatanabe} by 
\begin{equation*}
\left\langle M^t(\varphi) \right\rangle_T =  t^{3/2}\nu_t \int_{[0,T]\times \mathbb{R}} f(S,Y,\om)^2 \, \rmd S \rmd Y.
\end{equation*}
Now, by Cauchy-Schwarz inequality, we have:
\[f(S,Y,\om)^2 \leq \frac{r_t \lambda^2}{2\sqrt t}\int_{\mathbb R} \mathcal{Y}_t(S-,X)^2 \varphi(X)^2 \mathbf{1}_{|X-Y|\leq r_t/2\sqrt t} \, \rmd X,\]
so for all $U\leq T$ in $[0,1]$ and integer $p\geq 2$, we get by Tonelli's theorem that
\begin{align*}
& \IP\left[ \left| \left\langle M^t(\varphi) \right\rangle_T - \left\langle M^t(\varphi) \right\rangle_U\right|^p \right] \\
& \leq \frac{t^{3p/2} \nu_t^p r_t^p \lambda^{2p}}{2^p t^{p/2}} \IP\left[\int_{[U,T]\times \mathbb R} \mathcal{Y}_t(S-,X)^2 \varphi(X)^2 \left( \int_\mathbb{R} \mathbf{1}_{|X-Y|\leq r_t/2\sqrt t}\, \rmd Y \right) \rmd S \rmd X \right]^p\\
& = 2^{-p}t^{p/2} \nu_t^p r_t^{2p} \lambda^{2p} \, \IP\left[\int_{[U,T]\times \mathbb R} \mathcal{Y}_t(S-,X)^2 \varphi(X)^2 \, \rmd S \rmd X \right]^p.
\end{align*}

Then, observe that by our scaling relations \eqref{eq:paramRI}, the quantity $t^{1/2} \nu_t r_t \lambda(\b_t)^2$ is bounded by some constant $C>0$. Expanding the power of the integral, we get:
\begin{equation*}
\IP\left[ \left| \left\langle M^t(\varphi) \right\rangle_T - \left\langle M^t(\varphi) \right\rangle_U\right|^p \right] \leq C  \Vert \varphi \Vert^{2p}_\infty \int_{[U,T]^p\times\mathbb{R}^p}\IP\left[\prod_{i=1}^p \mathcal{Y}_t(S_i-,X_i)^2\right] \rmd \mathbf{S}\rmd \mathbf{X}
\end{equation*}
As we know by Lemma \ref{lem:boundOnMomentsOfwt} that $\mathcal{Y}_t(S,X)/\rho(S,X)$ is bounded in $L^{2p}$, we can again use the generalized H\"older inequality to bound the expectation of the product, and obtain that there is constant $C'>0$ such that
\begin{align*}
\IP\left[ \left| \left\langle M^t(\varphi) \right\rangle_T - \left\langle M^t(\varphi) \right\rangle_U\right|^p \right] & \leq C'  \Vert \varphi \Vert^{2p}_\infty \int_{[U,T]^p\times\mathbb{R}^p}\prod_{i=1}^p \rho(S_i,X_i)^2 \rmd \mathbf{S}\rmd \mathbf{X}\\
& = C'  \Vert \varphi \Vert^{2p}_\infty \left( \int_{[U,T]\times\mathbb{R}} \rho(S,X)^2\rmd S \rmd X \right)^p\\
&= C'\Vert \varphi \Vert^{2p}_\infty \, \pi^{-p/2}  \, | T^{1/2} - U^{1/2} |^p.
\end{align*}
Thus, Kolmogorov's tightness criterion \cite[Theorem 12.3]{billingsley1968convergence} applies, so the bracket $\left\langle M^t(\varphi) \right\rangle$ is tight. Hence Aldous' criterion applies, which concludes the proof of tightness of $M^t(\varphi)$.

\paragraph{Acknowledgments}
I would like to thank my PhD supervisor Francis Comets for suggesting this question and for his careful reading and helpful comments.

\bibliography{polymeres-bib.bib}

\begin{thebibliography}{56}
\makeatletter
\newcommand{\dinatlabel}[1]%
{\ifNAT@numbers\else\NAT@biblabelnum{#1}\hspace{2\labelsep}\fi}
\makeatother
\expandafter\ifx\csname natexlab\endcsname\relax\def\natexlab#1{#1}\fi
\expandafter\ifx\csname url\endcsname\relax\def\url#1{\texttt{#1}}\fi

\bibitem[Alberts u.\,a.(2014{\natexlab{a}})Alberts, Khanin und
  Quastel]{AlbertsKhaninQuastelCRP}
\dinatlabel{Alberts u.\,a. 2014{\natexlab{a}}} \textsc{Alberts}, Tom~;
  \textsc{Khanin}, Konstantin~; \textsc{Quastel}, Jeremy:
\newblock The continuum directed random polymer.
\newblock In: \emph{J. Stat. Phys.}
\newblock 154 (2014), Nr.~1-2, S.~305--326. --
\newblock URL \url{http://dx.doi.org/10.1007/s10955-013-0872-z}. --
\newblock ISSN 0022-4715

\bibitem[Alberts u.\,a.(2014{\natexlab{b}})Alberts, Khanin und
  Quastel]{AlbertsKhaninQuastelAP}
\dinatlabel{Alberts u.\,a. 2014{\natexlab{b}}} \textsc{Alberts}, Tom~;
  \textsc{Khanin}, Konstantin~; \textsc{Quastel}, Jeremy:
\newblock The intermediate disorder regime for directed polymers in dimension
  {$1+1$}.
\newblock In: \emph{Ann. Probab.}
\newblock 42 (2014), Nr.~3, S.~1212--1256. --
\newblock URL \url{http://dx.doi.org/10.1214/13-AOP858}. --
\newblock ISSN 0091-1798

\bibitem[Amir u.\,a.(2011)Amir, Corwin und Quastel]{ACQ11}
\dinatlabel{Amir u.\,a. 2011} \textsc{Amir}, Gideon~; \textsc{Corwin}, Ivan~;
  \textsc{Quastel}, Jeremy:
\newblock Probability distribution of the free energy of the continuum directed
  random polymer in {$1+1$} dimensions.
\newblock In: \emph{Comm. Pure Appl. Math.}
\newblock 64 (2011), Nr.~4, S.~466--537. --
\newblock URL \url{http://dx.doi.org/10.1002/cpa.20347}. --
\newblock ISSN 0010-3640

\bibitem[Bates und Chatterjee(2016)]{BatesChatterjee}
\dinatlabel{Bates und Chatterjee 2016} \textsc{Bates}, Erik~;
  \textsc{Chatterjee}, Sourav:
\newblock The endpoint distribution of directed polymers.
\newblock In: \emph{arXiv}
\newblock (2016)

\bibitem[Bertini und Cancrini(1995)]{BertiniCancrini95}
\dinatlabel{Bertini und Cancrini 1995} \textsc{Bertini}, Lorenzo~;
  \textsc{Cancrini}, Nicoletta:
\newblock The stochastic heat equation: {F}eynman-{K}ac formula and
  intermittence.
\newblock In: \emph{J. Statist. Phys.}
\newblock 78 (1995), Nr.~5-6, S.~1377--1401. --
\newblock URL \url{http://dx.doi.org/10.1007/BF02180136}. --
\newblock ISSN 0022-4715

\bibitem[Bertini und Giacomin(1997)]{BertiniGiacomin97}
\dinatlabel{Bertini und Giacomin 1997} \textsc{Bertini}, Lorenzo~;
  \textsc{Giacomin}, Giambattista:
\newblock Stochastic {B}urgers and {KPZ} equations from particle systems.
\newblock In: \emph{Comm. Math. Phys.}
\newblock 183 (1997), Nr.~3, S.~571--607. --
\newblock URL \url{http://dx.doi.org/10.1007/s002200050044}. --
\newblock ISSN 0010-3616

\bibitem[Billingsley(1968)]{billingsley1968convergence}
\dinatlabel{Billingsley 1968} \textsc{Billingsley}, Patrick:
\newblock \emph{Convergence of probability measures}.
\newblock John Wiley \& Sons, Inc., New York-London-Sydney, 1968. --
\newblock xii+253~S

\bibitem[Billingsley(1999)]{billing}
\dinatlabel{Billingsley 1999} \textsc{Billingsley}, Patrick:
\newblock \emph{Convergence of probability measures}.
\newblock Second.
\newblock John Wiley \& Sons, Inc., New York, 1999
\newblock (Wiley Series in Probability and Statistics: Probability and
  Statistics). --
\newblock x+277~S. --
\newblock URL \url{https://doi.org/10.1002/9780470316962}. --
\newblock A Wiley-Interscience Publication. --
\newblock ISBN 0-471-19745-9

\bibitem[Borodin u.\,a.(2014)Borodin, Corwin und Ferrari]{BoCorFe14}
\dinatlabel{Borodin u.\,a. 2014} \textsc{Borodin}, Alexei~; \textsc{Corwin},
  Ivan~; \textsc{Ferrari}, Patrik:
\newblock Free energy fluctuations for directed polymers in random media in
  {$1+1$} dimension.
\newblock In: \emph{Comm. Pure Appl. Math.}
\newblock 67 (2014), Nr.~7, S.~1129--1214. --
\newblock URL \url{https://doi.org/10.1002/cpa.21520}. --
\newblock ISSN 0010-3640

\bibitem[Caravenna u.\,a.(2017)Caravenna, Sun und Zygouras]{PolynomialChaos17}
\dinatlabel{Caravenna u.\,a. 2017} \textsc{Caravenna}, Francesco~;
  \textsc{Sun}, Rongfeng~; \textsc{Zygouras}, Nikos:
\newblock Polynomial chaos and scaling limits of disordered systems.
\newblock In: \emph{J. Eur. Math. Soc. (JEMS)}
\newblock 19 (2017), Nr.~1, S.~1--65. --
\newblock URL \url{https://doi.org/10.4171/JEMS/660}. --
\newblock ISSN 1435-9855

\bibitem[Comets(2017)]{CStFlour}
\dinatlabel{Comets 2017} \textsc{Comets}, Francis:
\newblock \emph{Lecture Notes in Mathematics}. Bd. 2175: \emph{Directed
  polymers in random environments}.
\newblock Springer, Cham, 2017. --
\newblock xv+199~S. --
\newblock URL \url{https://doi.org/10.1007/978-3-319-50487-2}. --
\newblock Lecture notes from the 46th Probability Summer School held in
  Saint-Flour, 2016. --
\newblock ISBN 978-3-319-50486-5; 978-3-319-50487-2

\bibitem[Comets und Yoshida(2004)]{CYkokyuroku}
\dinatlabel{Comets und Yoshida 2004} \textsc{Comets}, Francis~;
  \textsc{Yoshida}, Nobuo:
\newblock Some new results on {B}rownian directed polymers in random
  environment.
\newblock In: \emph{RIMS Kokyuroku}
\newblock 1386 (2004), S.~50--66

\bibitem[Comets und Yoshida(2005)]{CY05}
\dinatlabel{Comets und Yoshida 2005} \textsc{Comets}, Francis~;
  \textsc{Yoshida}, Nobuo:
\newblock Brownian directed polymers in random environment.
\newblock In: \emph{Comm. Math. Phys.}
\newblock 254 (2005), Nr.~2, S.~257--287. --
\newblock URL \url{http://dx.doi.org/10.1007/s00220-004-1203-7}. --
\newblock ISSN 0010-3616

\bibitem[Comets und Yoshida(2006)]{CY06}
\dinatlabel{Comets und Yoshida 2006} \textsc{Comets}, Francis~;
  \textsc{Yoshida}, Nobuo:
\newblock Directed polymers in random environment are diffusive at weak
  disorder.
\newblock In: \emph{Ann. Probab.}
\newblock 34 (2006), Nr.~5, S.~1746--1770. --
\newblock URL \url{https://doi.org/10.1214/009117905000000828}. --
\newblock ISSN 0091-1798

\bibitem[Comets und Yoshida(2013)]{CYBMPO2}
\dinatlabel{Comets und Yoshida 2013} \textsc{Comets}, Francis~;
  \textsc{Yoshida}, Nobuo:
\newblock Localization transition for polymers in {P}oissonian medium.
\newblock In: \emph{Comm. Math. Phys.}
\newblock 323 (2013), Nr.~1, S.~417--447. --
\newblock URL \url{http://dx.doi.org/10.1007/s00220-013-1744-8}. --
\newblock ISSN 0010-3616

\bibitem[Corwin(2016)]{corwin2016kardar}
\dinatlabel{Corwin 2016} \textsc{Corwin}, Ivan:
\newblock Kardar-{P}arisi-{Z}hang universality.
\newblock In: \emph{Notices of the AMS}
\newblock 63 (2016), Nr.~3, S.~230--239

\bibitem[Corwin und Tsai(2017)]{corwin2017}
\dinatlabel{Corwin und Tsai 2017} \textsc{Corwin}, Ivan~; \textsc{Tsai},
  Li-Cheng:
\newblock {KPZ} equation limit of higher-spin exclusion processes.
\newblock In: \emph{Ann. Probab.}
\newblock 45 (2017), 05, Nr.~3, S.~1771--1798. --
\newblock URL \url{https://doi.org/10.1214/16-AOP1101}

\bibitem[Dembo und Tsai(2016)]{Dembo2016}
\dinatlabel{Dembo und Tsai 2016} \textsc{Dembo}, Amir~; \textsc{Tsai},
  Li-Cheng:
\newblock Weakly asymmetric non-simple exclusion process and the
  {K}ardar-{P}arisi-{Z}hang equation.
\newblock In: \emph{Comm. Math. Phys.}
\newblock 341 (2016), Nr.~1, S.~219--261. --
\newblock URL \url{https://doi.org/10.1007/s00220-015-2527-1}. --
\newblock ISSN 0010-3616

\bibitem[Diehl u.\,a.(2017)Diehl, Gubinelli und Perkowski]{Diehl2017}
\dinatlabel{Diehl u.\,a. 2017} \textsc{Diehl}, Joscha~; \textsc{Gubinelli},
  Massimiliano~; \textsc{Perkowski}, Nicolas:
\newblock The {K}ardar-{P}arisi-{Z}hang equation as scaling limit of weakly
  asymmetric interacting {B}rownian motions.
\newblock In: \emph{Comm. Math. Phys.}
\newblock 354 (2017), Nr.~2, S.~549--589. --
\newblock URL \url{https://doi.org/10.1007/s00220-017-2918-6}. --
\newblock ISSN 0010-3616

\bibitem[Franco u.\,a.(2016)Franco, Gon\c{c}alves und Simon]{Franco2016}
\dinatlabel{Franco u.\,a. 2016} \textsc{Franco}, Tertuliano~;
  \textsc{Gon\c{c}alves}, Patr\'\i~c.~; \textsc{Simon}, Marielle:
\newblock Crossover to the stochastic {B}urgers equation for the {WASEP} with a
  slow bond.
\newblock In: \emph{Comm. Math. Phys.}
\newblock 346 (2016), Nr.~3, S.~801--838. --
\newblock URL \url{https://doi.org/10.1007/s00220-016-2607-x}. --
\newblock ISSN 0010-3616

\bibitem[Gon{\c{c}}alves und Jara(2014)]{Goncalves2014}
\dinatlabel{Gon{\c{c}}alves und Jara 2014} \textsc{Gon{\c{c}}alves},
  Patr{\'i}cia~; \textsc{Jara}, Milton:
\newblock Nonlinear Fluctuations of Weakly Asymmetric Interacting Particle
  Systems.
\newblock In: \emph{Archive for Rational Mechanics and Analysis}
\newblock 212 (2014), May, Nr.~2, S.~597--644. --
\newblock URL \url{https://doi.org/10.1007/s00205-013-0693-x}. --
\newblock ISSN 1432-0673

\bibitem[Gon\c{c}alves u.\,a.(2017)Gon\c{c}alves, Jara und
  Simon]{Goncalves2017}
\dinatlabel{Gon\c{c}alves u.\,a. 2017} \textsc{Gon\c{c}alves}, Patr\'icia~;
  \textsc{Jara}, Milton~; \textsc{Simon}, Marielle:
\newblock Second order {B}oltzmann-{G}ibbs principle for polynomial functions
  and applications.
\newblock In: \emph{J. Stat. Phys.}
\newblock 166 (2017), Nr.~1, S.~90--113. --
\newblock URL \url{https://doi.org/10.1007/s10955-016-1686-6}. --
\newblock ISSN 0022-4715

\bibitem[Gonçalves u.\,a.(2015)Gonçalves, Jara und Sethuraman]{goncalves2015}
\dinatlabel{Gonçalves u.\,a. 2015} \textsc{Gonçalves}, Patrícia~;
  \textsc{Jara}, Milton~; \textsc{Sethuraman}, Sunder:
\newblock A stochastic {B}urgers equation from a class of microscopic
  interactions.
\newblock In: \emph{Ann. Probab.}
\newblock 43 (2015), 02, Nr.~1, S.~286--338. --
\newblock URL \url{https://doi.org/10.1214/13-AOP878}

\bibitem[Gubinelli und Jara(2013)]{Gubinelli2013}
\dinatlabel{Gubinelli und Jara 2013} \textsc{Gubinelli}, M.~; \textsc{Jara},
  M.:
\newblock Regularization by noise and stochastic {B}urgers equations.
\newblock In: \emph{Stochastic Partial Differential Equations: Analysis and
  Computations}
\newblock 1 (2013), Jun, Nr.~2, S.~325--350. --
\newblock URL \url{https://doi.org/10.1007/s40072-013-0011-5}. --
\newblock ISSN 2194-041X

\bibitem[Gubinelli und Perkowski(2017)]{Gubinelli2017}
\dinatlabel{Gubinelli und Perkowski 2017} \textsc{Gubinelli}, Massimiliano~;
  \textsc{Perkowski}, Nicolas:
\newblock K{PZ} reloaded.
\newblock In: \emph{Comm. Math. Phys.}
\newblock 349 (2017), Nr.~1, S.~165--269. --
\newblock URL \url{https://doi.org/10.1007/s00220-016-2788-3}. --
\newblock ISSN 0010-3616

\bibitem[Gubinelli und Perkowski(2018)]{GubinelliEnergysolutions}
\dinatlabel{Gubinelli und Perkowski 2018} \textsc{Gubinelli}, Massimiliano~;
  \textsc{Perkowski}, Nicolas:
\newblock Energy solutions of {KPZ} are unique.
\newblock In: \emph{J. Amer. Math. Soc.}
\newblock 31 (2018), Nr.~2, S.~427--471. --
\newblock URL \url{https://doi.org/10.1090/jams/889}. --
\newblock ISSN 0894-0347

\bibitem[Hairer(2013)]{hairer}
\dinatlabel{Hairer 2013} \textsc{Hairer}, Martin:
\newblock Solving the {KPZ} equation.
\newblock In: \emph{Ann. of Math. (2)}
\newblock 178 (2013), Nr.~2, S.~559--664. --
\newblock URL \url{http://dx.doi.org/10.4007/annals.2013.178.2.4}. --
\newblock ISSN 0003-486X

\bibitem[Hairer und Shen(2017)]{hairer2017}
\dinatlabel{Hairer und Shen 2017} \textsc{Hairer}, Martin~; \textsc{Shen}, Hao:
\newblock A central limit theorem for the {KPZ} equation.
\newblock In: \emph{Ann. Probab.}
\newblock 45 (2017), 11, Nr.~6B, S.~4167--4221. --
\newblock URL \url{https://doi.org/10.1214/16-AOP1162}

\bibitem[Hoshino(2018)]{HOSHINO2017}
\dinatlabel{Hoshino 2018} \textsc{Hoshino}, Masato:
\newblock Paracontrolled calculus and {F}unaki--{Q}uastel approximation for the
  {KPZ} equation.
\newblock In: \emph{Stochastic Process. Appl.}
\newblock 128 (2018), Nr.~4, S.~1238--1293. --
\newblock URL \url{https://doi.org/10.1016/j.spa.2017.07.001}. --
\newblock ISSN 0304-4149

\bibitem[Huse und Henley(1985)]{HuHe85}
\dinatlabel{Huse und Henley 1985} \textsc{Huse}, David~A.~; \textsc{Henley},
  Christopher~L.:
\newblock Pinning and roughening of domain walls in {I}sing systems due to
  random impurities.
\newblock In: \emph{Physical review letters}
\newblock 54 (1985), Nr.~25, S.~2708

\bibitem[Ikeda und Watanabe(1989)]{IkedaWatanabe}
\dinatlabel{Ikeda und Watanabe 1989} \textsc{Ikeda}, Nobuyuki~;
  \textsc{Watanabe}, Shinzo:
\newblock \emph{North-Holland Mathematical Library}. Bd.~24: \emph{Stochastic
  differential equations and diffusion processes}.
\newblock Second.
\newblock North-Holland Publishing Co., Amsterdam; Kodansha, Ltd., Tokyo, 1989.
  --
\newblock xvi+555~S. --
\newblock ISBN 0-444-87378-3

\bibitem[Imbrie und Spencer(1988)]{ImbrieSpencer88}
\dinatlabel{Imbrie und Spencer 1988} \textsc{Imbrie}, John~; \textsc{Spencer},
  Thomas:
\newblock Diffusion of directed polymers in a random environment.
\newblock In: \emph{J. Statist. Phys.}
\newblock 52 (1988), Nr.~3-4, S.~609--626. --
\newblock URL \url{http://dx.doi.org/10.1007/BF01019720}. --
\newblock ISSN 0022-4715

\bibitem[Jacod und Shiryaev(2003)]{jacod2013limit}
\dinatlabel{Jacod und Shiryaev 2003} \textsc{Jacod}, Jean~; \textsc{Shiryaev},
  Albert~N.:
\newblock \emph{Grundlehren der Mathematischen Wissenschaften [Fundamental
  Principles of Mathematical Sciences]}. Bd. 288: \emph{Limit theorems for
  stochastic processes}.
\newblock Second.
\newblock Springer-Verlag, Berlin, 2003. --
\newblock xx+661~S. --
\newblock URL \url{https://doi.org/10.1007/978-3-662-05265-5}. --
\newblock ISBN 3-540-43932-3

\bibitem[Janson(1997)]{janson_1997}
\dinatlabel{Janson 1997} \textsc{Janson}, Svante:
\newblock \emph{Cambridge Tracts in Mathematics}. Bd. 129: \emph{Gaussian
  {H}ilbert spaces}.
\newblock Cambridge University Press, Cambridge, 1997. --
\newblock x+340~S. --
\newblock URL \url{https://doi.org/10.1017/CBO9780511526169}. --
\newblock ISBN 0-521-56128-0

\bibitem[Kardar u.\,a.(1986)Kardar, Parisi und Zhang]{KPZ}
\dinatlabel{Kardar u.\,a. 1986} \textsc{Kardar}, Mehran~; \textsc{Parisi},
  Giorgio~; \textsc{Zhang}, Yi-Cheng:
\newblock Dynamic Scaling of Growing Interfaces.
\newblock In: \emph{Physical Review Letters}
\newblock 56 (1986), Nr.~9, S.~889--892

\bibitem[Labb\'e(2017)]{Labbe2017}
\dinatlabel{Labb\'e 2017} \textsc{Labb\'e}, Cyril:
\newblock Weakly asymmetric bridges and the {KPZ} equation.
\newblock In: \emph{Comm. Math. Phys.}
\newblock 353 (2017), Nr.~3, S.~1261--1298. --
\newblock URL \url{https://doi.org/10.1007/s00220-017-2875-0}. --
\newblock ISSN 0010-3616

\bibitem[Last(2016)]{Last2016}
\dinatlabel{Last 2016} \textsc{Last}, G\"unter:
\newblock Stochastic analysis for {P}oisson processes.
\newblock In: \emph{Stochastic analysis for {P}oisson point processes} Bd.~7.
\newblock Bocconi Univ. Press, [place of publication not identified], 2016,
  S.~1--36. --
\newblock URL \url{https://doi.org/10.1007/978-3-319-05233-5_1}

\bibitem[Last und Penrose(2017)]{LastPenrose}
\dinatlabel{Last und Penrose 2017} \textsc{Last}, G\"unter~; \textsc{Penrose},
  Mathew:
\newblock \emph{Lectures on the {P}oisson process.}
\newblock Cambridge University Press (IMS Textbook), 2017

\bibitem[Mejane(2004)]{Mejane}
\dinatlabel{Mejane 2004} \textsc{Mejane}, Olivier:
\newblock Upper bound of a volume exponent for directed polymers in a random
  environment.
\newblock In: \emph{Ann. Inst. H. Poincar{\'e} Probab. Statist.}
\newblock 40 (2004), Nr.~3, S.~299--308. --
\newblock URL \url{http://dx.doi.org/10.1016/S0246-0203(03)00072-4}. --
\newblock ISSN 0246-0203

\bibitem[Mitoma(1983)]{mitoma1983tightness}
\dinatlabel{Mitoma 1983} \textsc{Mitoma}, Itaru:
\newblock Tightness of probabilities on {$C([0,1];{\cal S}^{\prime} )$}\ and
  {$D([0,1];{\cal S}^{\prime} )$}.
\newblock In: \emph{Ann. Probab.}
\newblock 11 (1983), Nr.~4, S.~989--999. --
\newblock URL
  \url{http://links.jstor.org/sici?sici=0091-1798(198311)11:4<989:TOPOA>2.0.CO;2-P&origin=MSN}.
  --
\newblock ISSN 0091-1798

\bibitem[Moreno~Flores(2014)]{Moreno14}
\dinatlabel{Moreno~Flores 2014} \textsc{Moreno~Flores}, Gregorio~R.:
\newblock On the (strict) positivity of solutions of the stochastic heat
  equation.
\newblock In: \emph{Ann. Probab.}
\newblock 42 (2014), Nr.~4, S.~1635--1643. --
\newblock URL \url{http://dx.doi.org/10.1214/14-AOP911}. --
\newblock ISSN 0091-1798

\bibitem[Mueller(1991)]{Mueller}
\dinatlabel{Mueller 1991} \textsc{Mueller}, Carl:
\newblock On the support of solutions to the heat equation with noise.
\newblock In: \emph{Stochastics Stochastics Rep.}
\newblock 37 (1991), Nr.~4, S.~225--245. --
\newblock URL \url{https://doi.org/10.1080/17442509108833738}. --
\newblock ISSN 1045-1129

\bibitem[Mukherjee u.\,a.(2016)Mukherjee, Shamov und
  Zeitouni]{Chiranjib-wkstrgdisorder}
\dinatlabel{Mukherjee u.\,a. 2016} \textsc{Mukherjee}, Chiranjib~;
  \textsc{Shamov}, Alexander~; \textsc{Zeitouni}, Ofer:
\newblock Weak and strong disorder for the stochastic heat equation and
  continuous directed polymers in {$d\geq 3$}.
\newblock In: \emph{Electron. Commun. Probab.}
\newblock 21 (2016), S.~Paper No. 61, 12. --
\newblock URL \url{https://doi.org/10.1214/16-ECP18}. --
\newblock ISSN 1083-589X

\bibitem[O'Connell und Yor(2001)]{OConnellYor}
\dinatlabel{O'Connell und Yor 2001} \textsc{O'Connell}, Neil~; \textsc{Yor},
  Marc:
\newblock Brownian analogues of {B}urke's theorem.
\newblock In: \emph{Stochastic Process. Appl.}
\newblock 96 (2001), Nr.~2, S.~285--304. --
\newblock URL \url{http://dx.doi.org/10.1016/S0304-4149(01)00119-3}. --
\newblock ISSN 0304-4149

\bibitem[Petermann(2000)]{Petermann}
\dinatlabel{Petermann 2000} \textsc{Petermann}, Markus:
\newblock \emph{Superdiffusivity of polymers in random environment}, Univ.
  Z\"urich, Dissertation, 2000

\bibitem[Sepp\"al\"ainen(2012)]{Sepp12}
\dinatlabel{Sepp\"al\"ainen 2012} \textsc{Sepp\"al\"ainen}, Timo:
\newblock Scaling for a one-dimensional directed polymer with boundary
  conditions.
\newblock In: \emph{Ann. Probab.}
\newblock 40 (2012), S.~19--73

\bibitem[Shiozawa(2009{\natexlab{a}})]{Shiozawa-clt}
\dinatlabel{Shiozawa 2009{\natexlab{a}}} \textsc{Shiozawa}, Yuichi:
\newblock Central limit theorem for branching {B}rownian motions in random
  environment.
\newblock In: \emph{J. Stat. Phys.}
\newblock 136 (2009), Nr.~1, S.~145--163. --
\newblock URL \url{http://dx.doi.org/10.1007/s10955-009-9774-5}. --
\newblock ISSN 0022-4715

\bibitem[Shiozawa(2009{\natexlab{b}})]{Shiozawa-loc}
\dinatlabel{Shiozawa 2009{\natexlab{b}}} \textsc{Shiozawa}, Yuichi:
\newblock Localization for branching {B}rownian motions in random environment.
\newblock In: \emph{Tohoku Math. J. (2)}
\newblock 61 (2009), Nr.~4, S.~483--497. --
\newblock URL \url{http://dx.doi.org/10.2748/tmj/1264084496}. --
\newblock ISSN 0040-8735

\bibitem[Tracy und Widom(1994)]{tracy1994level}
\dinatlabel{Tracy und Widom 1994} \textsc{Tracy}, Craig~A.~; \textsc{Widom},
  Harold:
\newblock Level-spacing distributions and the {A}iry kernel.
\newblock In: \emph{Comm. Math. Phys.}
\newblock 159 (1994), Nr.~1, S.~151--174. --
\newblock URL \url{http://projecteuclid.org/euclid.cmp/1104254495}. --
\newblock ISSN 0010-3616

\bibitem[Tracy und Widom(2008)]{Tracy2008}
\dinatlabel{Tracy und Widom 2008} \textsc{Tracy}, Craig~A.~; \textsc{Widom},
  Harold:
\newblock A {F}redholm determinant representation in {ASEP}.
\newblock In: \emph{J. Stat. Phys.}
\newblock 132 (2008), Nr.~2, S.~291--300. --
\newblock URL \url{https://doi.org/10.1007/s10955-008-9562-7}. --
\newblock ISSN 0022-4715

\bibitem[Tracy und Widom(2009)]{Tracy2009}
\dinatlabel{Tracy und Widom 2009} \textsc{Tracy}, Craig~A.~; \textsc{Widom},
  Harold:
\newblock Asymptotics in {ASEP} with step initial condition.
\newblock In: \emph{Comm. Math. Phys.}
\newblock 290 (2009), Nr.~1, S.~129--154. --
\newblock URL \url{https://doi.org/10.1007/s00220-009-0761-0}. --
\newblock ISSN 0010-3616

\bibitem[Tracy und Widom(2011)]{Tracy2008Erra}
\dinatlabel{Tracy und Widom 2011} \textsc{Tracy}, Craig~A.~; \textsc{Widom},
  Harold:
\newblock Erratum to: {I}ntegral formulas for the asymmetric simple exclusion
  process [MR2386729].
\newblock In: \emph{Comm. Math. Phys.}
\newblock 304 (2011), Nr.~3, S.~875--878. --
\newblock URL \url{https://doi.org/10.1007/s00220-011-1249-2}. --
\newblock ISSN 0010-3616

\bibitem[Vargas(2006)]{Va04}
\dinatlabel{Vargas 2006} \textsc{Vargas}, Vincent:
\newblock A local limit theorem for directed polymers in random media: the
  continuous and the discrete case.
\newblock In: \emph{Ann. Inst. H. Poincar{\'e} Probab. Statist.}
\newblock 42 (2006), Nr.~5, S.~521--534. --
\newblock URL \url{http://dx.doi.org/10.1016/j.anihpb.2005.08.002}. --
\newblock ISSN 0246-0203

\bibitem[Walsh(1986)]{Walsh}
\dinatlabel{Walsh 1986} \textsc{Walsh}, John~B.:
\newblock An introduction to stochastic partial differential equations.
\newblock In: \emph{\'{E}cole d'{\'e}t{\'e} de probabilit{\'e}s de
  {S}aint-{F}lour, {XIV}---1984} Bd.~1180.
\newblock Springer, Berlin, 1986, S.~265--439. --
\newblock URL \url{http://dx.doi.org/10.1007/BFb0074920}

\bibitem[W\"uthrich(1998)]{Wuthrich-fluctuation}
\dinatlabel{W\"uthrich 1998} \textsc{W\"uthrich}, Mario~V.:
\newblock Fluctuation results for {B}rownian motion in a {P}oissonian
  potential.
\newblock In: \emph{Ann. Inst. H. Poincar\'e Probab. Statist.}
\newblock 34 (1998), Nr.~3, S.~279--308. --
\newblock URL \url{https://doi.org/10.1016/S0246-0203(98)80013-7}. --
\newblock ISSN 0246-0203

\bibitem[W{{\"u}}thrich(1998)]{Wuthrich1}
\dinatlabel{W{{\"u}}thrich 1998} \textsc{W{{\"u}}thrich}, Mario~V.:
\newblock Superdiffusive behavior of two-dimensional {B}rownian motion in a
  {P}oissonian potential.
\newblock In: \emph{Ann. Probab.}
\newblock 26 (1998), Nr.~3, S.~1000--1015. --
\newblock URL \url{http://dx.doi.org/10.1214/aop/1022855742}. --
\newblock ISSN 0091-1798

\end{thebibliography}
\end{document}